\documentclass[AMA,Times1COL]{WileyNJDv5}

\articletype{Research Article}%

\received{Date Month Year}
\revised{Date Month Year}
\accepted{Date Month Year}
\journal{International Journal Of Robust And Nonlinear Control}
\volume{00}
\copyyear{2023}
\startpage{1}
\usepackage{graphicx}
 \usepackage{epsfig}
  \usepackage{lmodern}
\allowdisplaybreaks[4]
\begin{document}

\title{Block Backstepping Control for Rapid Stabilization of Multilayer Timoshenko Beams}

\author[1]{Guangwei Chen}
\author[2]{Rafael Vazquez}
\author[1]{Junfei Qiao}
\author[3]{Miroslav Krstic}

\authormark{Chen \textsc{et al.}}
\titlemark{Block Backstepping Control for Rapid Stabilization of Multilayer Timoshenko Beams}

\address[1]{\orgdiv{Beijing Laboratory of Smart Environmental Protection, Beijing Institute of Artificial Intelligence, Faculty of Information Technology}, \orgname{Beijing University of Technology}, \orgaddress{\state{Beijing}, \country{China}}}
\address[2]{\orgdiv{Dept. of Aerospace Engineering}, \orgname{ Universidad de Sevilla}, \orgaddress{\state{Sevilla}, \country{Spain}}}
\address[3]{\orgdiv{Dept. of Mechanical and Aerospace Engineering}, \orgname{University of California San Diego}, \orgaddress{\state{La Jolla}, \country{USA}}}
\corres{Corresponding author Guangwei Chen. \email{guangweichen@bjut.edu.cn}}

\presentaddress{100 Pingle Yuan, Chaoyang District, 100124, Beijing, China.}


\abstract[Abstract]{{In this paper, we consider the rapid stabilization of $N$-layer Timoshenko composite beams. Similar to our previous work on a two-layer composite beam, a Riemann-type transformation is proposed and used to transform the system with any number of layers into a 1-D hyperbolic PIDE-ODE form. However, when some of the layers have the same physical properties, the approach leads to isotachic hyperbolic PDEs which brings obstacles to controller design. After formulating the multilayer Timoshenko beam stabilization problem, this work considers the more general case of having a system with $(m+n)$ hyperbolic PIDEs and $\nu$ ODEs, and extends the theory of backstepping control of this general case to blocks of isotachic states, leading to a block backstepping design. Applying the result to multilayer Timoshenko beams, a boundary controller is obtained, achieving closed-loop stability of the origin in the $L^2$ sense. An arbitrarily rapid convergence rate can be obtained by adjusting control parameters. Finally, numerical simulations are presented corroborating the theoretical developments.}}


\maketitle
\footnotetext{\textbf{Abbreviations:} ODE, Ordinary Differential Equation; PIDE, Partial integral differential equation.}

\renewcommand\thefootnote{\fnsymbol{footnote}}
\setcounter{footnote}{1}

\section{Introduction}

\paragraph*{Composite Beams}
Composite beams have been widely used in various fields, such as  aeronautics~\cite{Ferreira2003}, mechanism design \cite{Attwood2016}, civil engineering~\cite{Barbero1991} or electronics~\cite{Waisman2002}. Several reasons justify their application, including weight reduction, higher overall stiffness,  enhanced properties (with respect to fracture, fatigue or corrosion) or cost reduction. However, when considering multi-layer beams, the coupling between layers can lead to vibration problems and, more critically, tip boundary conditions with anti-damping or anti-stiffness can cause divergence of the displacements and a consequent delamination of the beam into its unbonded constitutive layers. Therefore, it is necessary in some cases to design feedback controllers able to stabilize the equilibrium of the system.

\paragraph*{Literature: Stabilization of Timoshenko and Laminated Beams}
Regarding the single-layer Timoshenko beam, there exists plenty of literature on stabilization. For instance, Mattioni et al~\cite{mattioni2023lyapunov} propose a Lyapunov function with crossing terms in order to prove the exponential stability of a Timoshenko beam in the case of spatially varying
parameters and viscous damping in both the vertical and rotational
dynamics with the help of the port-Hamiltonian framework. Zhao et al~\cite{zhao2020finite}focus on a new finite-time convergence
disturbance rejection control scheme design for a flexible
Timoshenko manipulator subject to extraneous disturbances by developing  three disturbance
observers (DOs) and boundary controllers. Guo and Meng~\cite{guo2021robust} consider two-dimensional robust output tracking for a Timoshenko beam equation where the disturbances and references are produced by
a finite-dimensional exosystem, in which the observer-based error feedback control approach is adopted to achieve bounded control. Homaeinezhad and Abbasi~\cite{homaeinezhad2023feedback} deal with tracking control and active vibration suppression of a
hybrid rigid-elastic vibrational mobile system with a nonlinear elastic foundation in the presence of actuator bandwidth limitation and control input saturation. Meng et al~\cite{meng2023reduced} propose observer based controls for a MIMO Timoshenko beam system with disturbances and references from an exosystem, and by further limiting the references and disturbances, also gives reduced order controls. Wang et al~\cite{wang2023saturated} consider the saturated feedback
stabilization for a nonuniform Timoshenko beam with boundary
control matched unknown internal uncertainties and external disturbances in which an infinite-dimensional disturbance estimator without involving high gain is designed to estimate the total disturbances. Redaud et al~\cite{redaud2022domain} combine the backstepping
methodology and the port Hamiltonian framework to design a boundary full-state feedback controller that modifies
the closed-loop in-domain damping of a Timoshenko beam.

Designing controllers for $N$-layer beams is more challenging due to the increased complexity of the underlying mathematical models, but there also exists some literature. For example, Lo and Tatar~\cite{Lo2015} studied the stabilization of a laminated beam with interfacial slip, with an adhesive of small thickness bonding the two layers and creating damping. In another work~\cite{Mustafa2018}, a viscoelastic laminated beam model is considered without additional control, and   explicit energy decay formulae are established,  giving the optimal decay rates by using minimal  conditions on a relaxation function. The works~\cite{Apalara2020, Apalara2020_1} investigated a one-dimensional laminated Timoshenko beam with a single nonlinear structural damping due to interfacial slip, and established an explicit and general decay result by adopting a multiplier method exploiting some properties of convex functions. Guesmia et al~\cite{Guesmia2021} researched the well-posedness and stability of  structures with interfacial slip; a large class of control kernels are considered and the system is proven to have a unique solution satisfying certain regularity properties. Apart from damping control, boundary control are frequently applied to laminated Timoshenko beams, for instance, Cao et al~\cite{Cao2007} considered the stability of the closed loop system composed of laminated beams with boundary feedback controls, and a simple test method was used to verify exponential stability. Some researchers have adopted the simultaneous use of  boundary  and interfacial damping control to obtain  exponential stability~\cite{Tatar2015, Alves2020}. Considering the time delay, Feng~\cite{Feng2018} studied the long-time dynamics of laminated Timoshenko beams and established the existence of smooth finite-dimensional global attractors for the corresponding solution semigroup.

\paragraph*{Literature: Rapid Stabilization of Hyperbolic PDEs}
Most of these works achieve stability or even exponential stability, but not \emph{rapid} stabilization (being able to set an arbitrarily fast decay rate in the closed loop), much less in the presence of destabilizing boundary conditions. This goal was achieved for beams by the use of backstepping technique in two pioneering works~\cite{krstic2008control,smyshlyaev2009arbitrary}, in the first case for an undamped shear beam, and in the second for an Euler-Bernoulli beam. More recently, backstepping was extended to  obtain rapid stabilization for one-layer Timoshenko~\cite{chen2022} and the two-layer Timoshenko~\cite{chen2023} (with potentially destabilizing boundary conditions) by using a Riemann transformation to cast the system as a 1-D hyperbolic PIDE-ODE system. This allows the use of the backstepping method, which has provided many designs for hyperbolic systems over the years; starting from a single 1-D hyperbolic partial integro-differential equation~\cite{krstic3}, the method was to extended next to $2 \times 2$ systems ~\cite{vazquez-linear,vazquez-nonlinear}, to $n+1 \times n+1$ system~\cite{florent}s ($n$ states convecting in one direction with one counter-convecting state that is controlled), and finally to the general case of $n+m \times n+m$ systems ($n$ states convecting in one direction and $m$ controlled states convecting in the opposite), both in the linear~\cite{hu} and quasi-linear~\cite{c6} cases. A later refinement allowed to obtain minimum-time convergence~\cite{Auriol2016, espitia2021predictor}. {
Block backstepping is developed for the case where some of the layers possibly have the same physical properties (as e.g. in lamination of repeated layers) which leads to \emph{isotachic} hyperbolic PDEs (i.e. where some states have the same transport speeds). The particular yet physical and interesting case of Timoshenko beams with repeated layers requires modification of the general design, which we provide as the traditional PDE backstepping suitable for hyperbolic PDEs with different transport speeds is not aplicable. This modifications allow to use backstepping PIDE-ODE designs~\cite{c1,deutscher2019output,gabriel2023robust}, which we applied for one a two-layer Timoshenko beam designs.\\ Regarding block backstepping, the literature is scarce; for example, Humaidi et al~\cite{humaidi2018design} proposed two block-backstepping algorithms for balancing and tracking control of ball and arc system and showed the block-backstepping controller designed for nonlinear system gives better transient performance. Another work~\cite{humaidi2017design} presented three block-backstepping controllers for ball and arc system and the internal stability of the system was analyzed using zero dynamic criteria. Similarly, for controlling the Two-Wheeled Inverted Pendulum (TWIP) system, Humaidi et al~\cite{humaidi2021block} studied the the control design of block backstepping control and the  simulation showed that the proposed controller achieves the  goals.}

\paragraph*{Contribution} The contribution and novelty of this paper is threefold:  
{
\begin{itemize}
    \item First, the  application of a Riemann-like transformation to transform the Timoshenko beam equations into a PIDE-ODE hyperbolic system (and then design a backstepping control law) is generalized from basic Timoshenko beams to Timoshenko beams with any number of layers (under spatially-varying coefficients), showing that the method we proposed has quite wide adaptability for multi-layer beams with any structural characteristics; 
    \item Second,  block backstepping is developed for the case where some of the layers possibly have the same physical properties (as e.g. in lamination of repeated layers) leading to \emph{isotachic} hyperbolic PDEs (i.e. where some states have the same transport speeds). The particular yet physical and interesting case of Timoshenko beams with repeated layers requires modification of the general design, which we provide. The traditional PDE backstepping design, suitable for hyperbolic PDEs with different transport speeds, does not work here. Therefore, the block backstepping needs to be developed, considering ``blocks'' of those PDE states having an identical transport speed. After reducing the blocks, PIDE-ODE designs~\cite{c1,deutscher2019output,gabriel2023robust} can be applied; 
    \item Third, rapid stabilization, with an arbitrary decay rate imposed by adjusting parameters in the control laws, is achieved for $N$-layer beams for the first time in this paper. The closed-loop exponential stability with an arbitrary rate of convergence is proved by Lyapunov-based analysis. Finally, the numerical solutions of a rather complex system with discontinuos kernel equations are obtained by power series method, and the numerical simulation verifies the theoretical result.
\end{itemize}}

The paper is organized as follows: Section \ref{sec2} presents the multi-layer Timoshenko beam and its representation as a 1-D $(m+n)$ hyperbolic PIDE system coupled with  $m$ ODEs,  involving some isotachic PDE states that are written as blocks. Section \ref{sec-control} gives the block design of the boundary controller for this case, and the main result. Then, Section \ref{contdesign} analyzes the resulting controller. Section \ref{stability} studies the closed-loop stability. Section \ref{sec-extension} applies the proposed method to the $N$-layer Timoshenko beams written in the PIDE-ODE form. Section \ref{num-simulation} shows a numerical simulation corroborating the theoretical results. Finally, Section \ref{sec-conclusion} closes the paper with some concluding remarks.


\section{Problem Statement}\label{sec2}
In this paper, we investigate the following $N$-layer Timoshenko composite beam model. {Before we introduce the model given in (\ref{plant_eq1})-(\ref{plant_C}) below, we point out that it is related to (\ref{plant_bc2}) in literature~\cite{Lenci2013}. However, the paper~\cite{Lenci2013} considers
only a two-layer beam and also takes the longitudinal displacements into
account. In contrast, (\ref{plant_eq1})-(\ref{plant_C}) involves more general couplings between the layers. Moreover, that paper ~\cite{Lenci2013} does not have boundary conditions as in (\ref{plant_bc2}).}

\begin{figure*}[t]
\begin{centering}
	 \includegraphics[width=13cm]{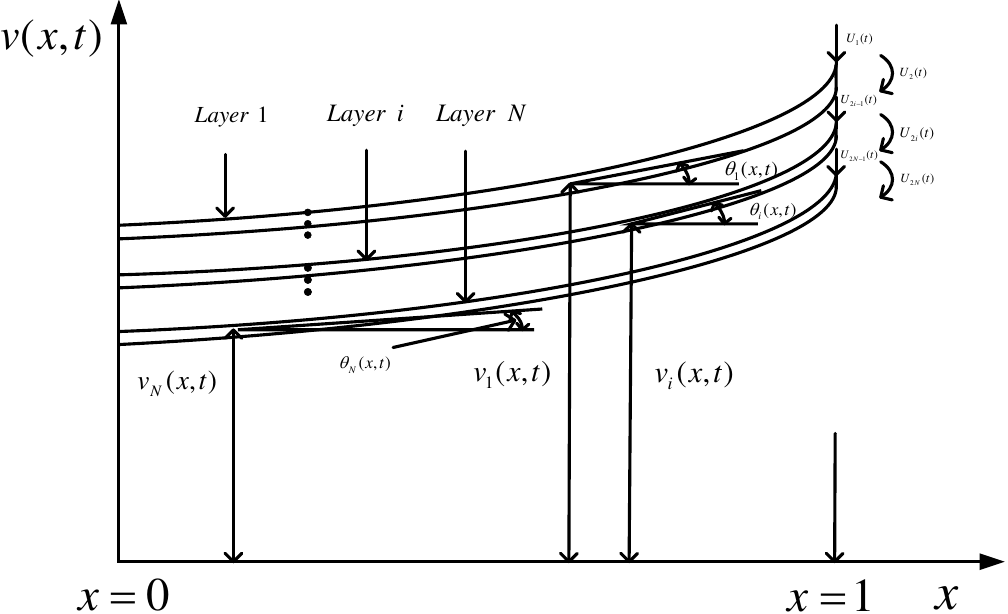}
		\caption{{Illustration of multilayer Timoshenko beam system}}
        \label{Fig-structure}
\end{centering}
\end{figure*}

So, we consider the $N$-layer Timoshenko beam model
\begin{align}
 \beta_i {v}_{i,tt}&=\eta_{i}\left(v_{i,xx}+\theta_{i,x}\right)+C_1(i)k^{i-1}_{no} s^{i-1}_{no}-C_2(i)k^i_{no} s^i_{no}\label{plant_eq1}\\
	\zeta_{i} {\theta}_{i,tt}&=\alpha_i\theta_{i,xx}-\eta_{i}\left(v_{i,x}+\theta_{i}\right)+C_1(i)h^{i-1}_{2} k^{i-1}_{ta} s^{i-1}_{ta}+C_2(i)h^i_{1} k^i_{ta} s^i_{ta},\label{plant_eq2}\\
	s^{i-1}_{ta}&=-h^{i-1}_{1} \theta_{i-1}-h^{i-1}_{2} \theta_{i}, \label{plant_eq3}\\
	s^{i-1}_{no}&= v_{i-1}-v_{i},~~i=1,\cdots,N\label{plant_eq4}
\end{align}
where 
\begin{align}
C_1(i)=\left\{
\begin{aligned}
    1 \quad i>1\\
    0 \quad i=1
    \end{aligned}
\right.~,C_2(i)=\left\{
\begin{aligned}
    1 \quad i<N\\
    0 \quad i=N
    \end{aligned}
\right.\label{plant_C}
\end{align}
where the sub-index $i$ makes reference to each of the layers, $v_i$ are  the transversal displacements, $\theta_i$ the  rotational angles of the cross-sections, $\eta_i$ the shear stiffnesses, $\zeta_i$ the rotational
inertia, $h_i$ the interface-centroids distances, $k^i_{ta}$ and $k^i_{no}$ the tangential and normal interface stiffnesses, $\alpha_i$ and $\beta_i$ the ratios of the $i$-layer beam with respect to its normal stiffnesses and its moments of inertia
of the cross-section, $s^i_{ta}$ and $s^i_{no}$ the tangential and normal displacements in the interface between beams, 
with boundary conditions
\begin{align}
	v_{i,x}(0,t)&=\theta_{i}(0,t)-\xi_{2i-1} v_{i,t}(0,t)-\xi_{2i} v_{i}(0,t), \notag\\
	v_{i,x}(1,t)&=U_{2i-1}(t), \notag\\
	\theta_{i,x}(0,t)&=0, \quad \theta_{i,x}(1,t)=U_{2i}(t),~~i=1,\cdots,N\label{plant_bc2}
\end{align}
where $\xi_{2i-1}$ are the anti-damping of each beam, and  $\xi_{2i}$  the anti-stiffness, with $U_1(t)$,   $U_2(t)$,  $U_3(t)$, $\cdots$, $U_{2N}(t)$ being the actuation variables that are forces and torques, which are designed in subsequent sections.  All quantities in the model are dimensionless. It must note that the $2N$ actuators are independent and they will not de-laminate the adjacent layers due to the adhesives existing between them. {Fig. \ref{Fig-structure} depicts the multilayer Timoshenko beam system.}

{Next, we consider a hyperbolic PIDE-ODE system capturing multilayer Timoshenko dynamics through a Riemann-like transformation. 
\begin{assumption} \label{assump2}
	The anti-damping coefficients $\xi_{2i-1}$ appearing in (\ref{plant_bc2}) verify $\xi_{2i-1} \neq \sqrt{\beta_i}/ \sqrt {\eta_{i}}, ~i=1,2,\cdots, N$.
\end{assumption}
This assumption, that only precludes very specific values of the coefficients, is required for the application of the Riemann transformation defined next; see the explanation in~\cite{chen2022} for the equivalent assumption in the 1-layer case.}

{
Following our previous work \cite{chen2022}, we introduce the transformations
\begin{align}\label{eq2}
 p_{i} &= \sqrt{\eta_i}v_{i,x}+\sqrt{\beta_i}v_{i,t},~r_{i}=\sqrt{\eta_i}v_{i,x}-\sqrt{\beta_i}v_{i,t},\\
	 q_{i}&=\sqrt{\alpha_i}\theta_{i,x}+\sqrt{\zeta_i}\theta_{i,t},~ s_{i}=\sqrt{\alpha_i}\theta_{i,x}-\sqrt{\zeta_i}\theta_{i,t}, ~i=1,2,\cdots,N.\label{Treq2}
\end{align}
and the definition of states
\begin{align}
    x_{2i-1}&=v_{i}(0,t),\\
 x_{2i}&=\theta_{i}(0,t), ~i=1,2,\cdots,N. \label{eq2-bis}
\end{align}
Then the $N$-layer Timoshenko composite beams are transformed into the PIDE-ODE system:
\begin{align}  \label{Equ_eq1}
	p_{i,t}=&\frac{\sqrt{\eta_i}}{\sqrt{\beta_i}}p_{i,x}+\frac{\eta_i}{2\sqrt{\alpha_i\beta_i}}(q_i+s_i)+C_1(i)\frac{k^{i-1}_n}{\sqrt{\beta_i}}\left(\int^x_0\frac{p_{i-1}+r_{i-1}}{2\sqrt{\eta_{i-1}}}-\frac{p_i+r_i}{2\sqrt{\eta_i}}ds+x_{2i-3}-x_{2i-1}\right)\notag\\
  &-C_2(i)\frac{k^i_n}{\sqrt{\beta_i}}\left(\int^x_0\frac{p_{i}+r_{i}}{2\sqrt{\eta_{i}}}-\frac{p_{i+1}+r_{i+1}}{2\sqrt{\eta_{i+1}}}ds+x_{2i-1}-x_{2i+1}\right),\\
	r_{i,t}=&-\frac{\sqrt{\eta_i}}{\sqrt{\beta_i}}r_{i,x}-\frac{\eta_i}{2\sqrt{\alpha_i\beta_i}}(q_i+s_i)-C_1(i)\frac{k^{i-1}_n}{\sqrt{\beta_i}}\left(\int^x_0\frac{p_{i-1}+r_{i-1}}{2\sqrt{\eta_{i-1}}}-\frac{p_i+r_i}{2\sqrt{\eta_i}}ds+x_{2i-3}-x_{2i-1}\right)\notag\\
  &+C_2(i)\frac{k^i_n}{\sqrt{\beta_i}}\left(\int^x_0\frac{p_{i}+r_{i}}{2\sqrt{\eta_{i}}}-\frac{p_{i+1}+r_{i+1}}{2\sqrt{\eta_{i+1}}}ds+x_{2i-1}-x_{2i+1}\right),\\
  q_{i,t}=&\frac{\sqrt{\alpha_i}}{\sqrt{\zeta_i}}q_{i,x}-\frac{\sqrt{\eta_i}}{2\sqrt{\zeta_i}}(p_i+r_i)-C_1(i)\frac{h^{i-1}_2k^{i-1}_th^{i-1}_1}{\sqrt{\zeta_i}}\left(\int^x_0 \frac{q_{i-1}+s_{i-1}}{2\sqrt{\eta_i}}dy+x_{2i-2}\right)\notag\\
	&-\left(\frac{\eta_i}{\sqrt{\zeta_i}}+C_1(i)\frac{(h^{i-1}_2)^2k^{i-1}_t}{\sqrt{\zeta_i}}+C_2(i)\frac{(h^{i}_1)^2k^{i}_t}{\sqrt{\zeta_i}}\right)\left(\int^x_0\frac{q_i+s_i}{2\sqrt{\eta_i}}ds+x_{2i}\right)\notag\\
	&-C_2(i)\frac{h^i_1h^i_2k^i_t}{\sqrt{\zeta_1}}\left(\int^x_0\frac{q_{i+1}+s_{i+1}}{2\sqrt{\eta_{i+1}}}ds+x_{2i+2}\right),\\ 
	 s_{i,t}=&-\frac{\sqrt{\alpha_i}}{\sqrt{\zeta_i}}s_{i,x}+
  \frac{\sqrt{\eta_i}}{2\sqrt{\zeta_i}}(p_i+r_i)+C_1(i)\frac{h^{i-1}_2k^{i-1}_th^{i-1}_1}{\sqrt{\zeta_i}}\left(\int^x_0 \frac{q_{i-1}+s_{i-1}}{2\sqrt{\eta_i}}dy+x_{2i-2}\right)\notag\\
	&+\left(\frac{\eta_i}{\sqrt{\zeta_i}}+C_1(i)\frac{(h^{i-1}_2)^2k^{i-1}_t}{\sqrt{\zeta_i}}+C_2(i)\frac{(h^{i}_1)^2k^{i}_t}{\sqrt{\zeta_i}}\right)\left(\int^x_0\frac{q_i+s_i}{2\sqrt{\eta_i}}ds+x_{2i}\right)\notag\\
	&+C_2(i)\frac{h^i_1h^i_2k^i_t}{\sqrt{\zeta_1}}\left(\int^x_0\frac{q_{i+1}+s_{i+1}}{2\sqrt{\eta_{i+1}}}ds+x_{2i+2}\right),\\ \label{ODE_eq1}
	\dot x_{2i-1}&=\frac{1}{\sqrt{\beta_i}-\xi_{2i-1}\sqrt{\eta_i}}(p_i(0)+\xi_{2i}\sqrt{\eta_i}x_{2i-1}-x_{2i}), \\
	\dot x_{2i}&=\frac{q_i(0)}{\sqrt{\zeta_i}}, i=1,2,\cdots,N,\label{ODE_eq2}
\end{align}
with the following boundary conditions:
\begin{align}\label{BD_eq8}
	r_i(0)&=\frac{\xi_{2i-1}\sqrt{\eta_i}+\sqrt{\beta_i}}{\xi_{2i-1}\sqrt{\eta_{i}}-\sqrt{\beta_i}}p_i(0)+\frac{2\sqrt{\beta_i}}{\sqrt{\beta_i}-\xi_{2i-1}\sqrt{\eta_{i}}}\left(x_{2i}-\xi_{2i}\sqrt{\eta_i}x_{2i-1}\right),\notag\\
	s_i(0)&=-q_i(0), \notag\\
	p_i(1)&=U_{i,p}, q_i(1)=U_{i,q}, i=1,2,\cdots, N.
\end{align}
where  ${U_{i,p}}(t)=\sqrt{\eta_i}U_{2i-1}(t)+\sqrt{\beta_i}{v_{i, t}}(1, t)$ and ${U_{i,q}(t)=\sqrt{\alpha_{i}}}U_{2i}(1,t)+\sqrt{\zeta_i}\theta_{i,t}(1,t)$ are the redefined control variables for this system.} 
{
Ordering the states $p_i$, $r_i$, $q_i$ and $s_i$ as function of their transport speeds, namely $\frac{\sqrt{\eta_i}}{\sqrt{\beta_i}},\frac{\sqrt{\alpha_i}}{\sqrt{\zeta_i}}, i=1,2,\dots,N$, is very important to implement the controller design. There is the possibility that some transport speeds are equal; in that case, the corresponding states can   be treated as a ``block''. We assume there are $\kappa$ different transport speeds $T_1>T_2>\cdots>T_{\kappa}$. According to the order of transport speeds $T_1, T_2, \cdots,T_{\kappa}$, we order the PDE states $p_i,r_i, q_i,s_i$ as $p^{or}_i,r^{or}_i,q^{or}_i,s^{or}_i, i=1,2,\cdots, N$ with the corresponding parameters $\eta^{or}_i,\beta^{or}_i,\alpha^{or}_i,\zeta^{or}_i$. Define
\begin{align}
	Y&=\left[
	\begin{array}{ccccc}
		p^{or}_1 &q^{or}_1 &\cdots&p^{or}_N&q^{or}_N\\
	\end{array}
	\right]^T, \\
	Z&=\left[
	\begin{array}{ccccc}
		r^{or}_1 &s^{or}_1 &\cdots&r^{or}_N&s^{or}_N\\
	\end{array}
	\right]^T,\\
	X&=\left[
	\begin{array}{ccccc}
		x_1 &x_2 &\cdots&x_{2N-1}&x_{2N}\\
	\end{array}
	\right]^T,\\
	U^{pq}&=\left[
	\begin{array}{ccccc}
		U^{or}_{1,p} &U^{or}_{1,q} &\cdots&U^{or}_{N,p}&U^{or}_{N,q}\\
	\end{array}
	\right]^T,
\end{align}
\begin{align}
	Z_t(x,t)&=-\Sigma Z_x(x,t) -\Lambda (Z(x,t)+Y(x,t))-\Pi X(t)-\int^x_0 F[Z(y,t)+Y(y,t)]dy,\label{Tm_e1}\\ 
	Y_t(x,t)&=\Sigma Y_x(x,t)+\Lambda (Y(x,t)+ Z(x,t))+\Pi X(t)+\int^x_0 F[Z(y,t)+Y(y,t)]dy,\label{Tm_e2}
	\\ 
	\dot{X}(t)&=AX(t)+BY(0, t),
\end{align}
with boundary conditions
\begin{align}
	Y(1,t)=U^{pq}(t),~	Z(0,t)=CY(0, t)+DX(t),\label{TBDm_e}
\end{align}
\begin{align} \label{eq4}
  \Sigma&=\mathrm{diag}\left\{\frac{\sqrt{\eta^{or}_1}}{\sqrt{\beta^{or}_1}},\frac{\sqrt{\alpha^{or}_1}}{\sqrt{\zeta^{or}_1}} ,\dots,\frac{\sqrt{\eta^{or}_N}}{\sqrt{\beta^{or}_N}},\frac{\sqrt{\alpha^{or}_N}}{\sqrt{\zeta^{or}_N}}\right\}\,,
 \end{align}
with the values in $\Sigma$ ordered from largest to smallest, and the other coefficient matrices $\Pi,F,\Lambda$ are similarly obtained.
In particular, the matrix $B$ consists of the coefficients of~$r_i(0,t)$, $s_i(0,t)$; since it is diagonal and invertible, it verifies Assumption \ref{assump1}. Inspired by the form of (\ref{Tm_e1})--(\ref{TBDm_e}), we start by solving the stabilization problem for a general structure of a system of $(m+n)$ hyperbolic PIDEs with $m$ ODEs constant-coefficient model, which is also an interesting problem by itself.} The plant allows the existence of isotachic states, which requires expanding the previous theory as explained next.

The system equations are
{\begin{align}
Z_t(x,t)&=-\Sigma^+(x) Z_x(x,t) +\Lambda^{++} (x)Z(x,t)+\Lambda^{+-}(x)Y(x,t)+\Pi^+(x)X(t)+\int^x_0 \left[F^{++}(y)Z(y,t)+F^{+-}(y)Y(y,t)\right]dy,\label{plant_eq1_Hyperbolic}\\ 
	Y_t(x,t)&=\Sigma^-(x) Y_x(x,t)+\Lambda^{-\,-}(x) Y(x,t)+ \Lambda^{-+}(x)Z(x,t)+\Pi^{-}(x) X(t)+\int^x_0 \left[F^{-+}(y)Z(y,t)+F^{-\,-}(y)Y(y,t)\right]dy,\label{plant_eq2_Hyperbolic}
	\\ 
	\dot{X}&=AX+BY(0, t),
 \end{align}}
with boundary conditions
\begin{align}
Y(1,t)=U(t),~Z(0,t)=CY(0, t)+DX(t). \label{BDm_e_hyperbolic}
\end{align}
where 
\begin{align}
    Z&=[z_1~z_2~\cdots~z_m]^T, Y=[y_1~y_2~\cdots~y_n]^T,X=[x_1~x_2~\cdots~x_{\nu}]^T\\
    \Sigma^{+}(x)&=\left[\begin{array}{ccc}
		\epsilon^+_1(x)&~&0\\
		~&\ddots&~\\
            0&~&~\epsilon^+_{m}(x)\\
	\end{array}
	\right]\in \mathbb{R}^{m\times m},\Sigma^{-}(x)=\left[\begin{array}{ccc}
		\Sigma^-_1(x)&~&0\\
		~&\ddots&~\\
            0&~&~\Sigma^-_{\kappa}(x)\\
	\end{array}
	\right]\in \mathbb{R}^{n\times n},\\
 \Sigma^-_j&=\epsilon^-_j(x)I_{n_j}, \sum_{j=1}^{\kappa} n_j=n,~n_j\in \mathbb{N}^+, \Pi^{+}(x), D\in\mathbb{R}^{m\times \nu}, \Lambda^{++}(x), F^{++}(x)\in \mathbb{R}^{m\times m}, A\in \mathbb{R}^{\nu\times\nu},\label{isotachic}\\
\Lambda^{+-}(x)&, F^{+-}(x)\in \mathbb{R}^{m\times n},B,C\in \mathbb{R}^{\nu\times n}, F^{-\,-}(x), \Lambda^{-\,-}(x)\in \mathbb{R}^{n\times n}, \Pi^{-}(x), \Lambda^{-+}(x), F^{-+}(x)\in \mathbb{R}^{n\times m},
\end{align}
with transport speeds
\begin{align}
-\epsilon^-_1(x)<\cdots<-\epsilon^-_{\kappa}(x)<0<\epsilon^+_1(x)\leq\cdots\leq\epsilon^+_m(x).
\end{align}

 In (\ref{isotachic}), $I_{n_j}$ represents the $n_j$-sized identity matrix. $\Sigma^-_j(x)$ are the blocks that assemble to $\Sigma^-(x)$. Thus, the $Y$-system has $\kappa$ different transport speeds. When $\kappa=n$, that is $n_1=\cdots=n_{\kappa}=1$, all the states of $Y-$system have different transport speeds (non-isotachic case) and the classical backstepping design~\cite{hu} (or its multiple variations) can be directly applied. When $\kappa<n$, there are at least two states having an identical transport speed (isotachic case). The direct coupling terms among the states with the same transport speed produce  singularities in the kernel equations~\cite{hu}. To guarantee the kernel equations are solvable, {we expand upon Remark 6 of~\cite{hu} (which was only considered for a hyperbolic PDE system without  integral terms and coupled ODEs) to this more general PIDE-ODE case that incorporates  spatially-varying coefficients as well} and introduce an invertible transformation $\mathcal A(x)$ to transform the original system into an intermediate system where the isotachic states have no  coupling between them. The details of transformation $\mathcal A(x)$ are presented in Section~\ref{sec-transformation}. 

In addition, the following assumption is essential to achieve the arbitrarily rapid stabilization of the coupled hyperbolic PIDE-ODE system.
\begin{assumption} \label{assump1}
    The matrix pair $(A,B)$ is controllable.
\end{assumption}


\section{Controller design}\label{sec-control}

\subsection{Block transformation for isotachic states}\label{sec-transformation}

For ease of derivation, the coupling matrix $\Lambda^{-\,-}(x)$ can be rewritten using blocks as follows
{
\begin{align}
\Lambda^{-\,-}(x)=\left[
	\begin{array}{ccc}
		\Lambda^{-\,-}_{11}(x)&\cdots&\Lambda^{-\,-}_{1\kappa}(x)\\
		\vdots&\ddots&\vdots\\
            \Lambda^{-\,-}_{\kappa 1}(x)& \cdots &~\Lambda^{-\,-}_{\kappa\kappa}(x)\\
	\end{array}
	\right]
\end{align}}
where $\Lambda^{-\,-}_i$ refers to the coupling only between states belonging to the $i$-th block of isotachic states.
To eliminate these, we introduce a transformation
\begin{align}\label{matrix_T1}
    \bar Y(x,t)&=\mathcal A(x){Y(x,t)}\\
     \mathcal A(x)&=\mathrm{diag} \{  \mathcal A_1(x),\mathcal A_2(x), \cdots,\mathcal A_\kappa(x) \}\\
     \frac{d\mathcal A_j(x)}{dx}&=\frac{1}{\sigma^{-}_{j}}\mathcal A_j(x){\Lambda^{-\,-}_{jj}(x)},~\mathcal A_j(0)=I_{n_j}
\end{align}
The matrices $\mathcal A_j(x)$ are all invertible, with their inverses $\tilde {\mathcal A}_j(x)=\left(\mathcal A_j(x)\right)^{-1}$ verifying
\begin{align}\label{matrix_T2}
    \frac{d\tilde {\mathcal A}_j(x)}{dx}=-\frac{1}{\sigma^-_j}{\Lambda^{-\,-}_{jj}(x)}\tilde{\mathcal A}_j(x),~\tilde{\mathcal A}_j(0)=I_{n_j}
\end{align}
It is easy to see  $\tilde {\mathcal A}_j(x)$ is the inverse transformation of ${\mathcal A}_j(x)$ since ${\mathcal A}_j(0)\tilde {\mathcal A}_j(0)=I_{n_j}$ and $\frac{d{\mathcal A}_j(x)\tilde {\mathcal A}_j(x)}{dx}=0$.
{
Use the transformation (\ref{matrix_T1})--(\ref{matrix_T2}) and define ${\bar\Lambda}^{+-}(x)=\Lambda^{+-}(x)\tilde{\mathcal {A}}(x), {\bar F}^{+-}(y)=F^{+-}(y)\tilde{\mathcal {A}}(y), \bar\Lambda^{-\,-}(x)={\mathcal {A}}(x)\left[\Lambda^{-\,-}(x)-\Sigma^{-}(x) \tilde{\mathcal{A}}(x) \frac{d \mathcal{A}}{d x}(x)\right] \tilde{\mathcal {A}}(x),$
$\bar \Lambda^{-+}(x)={\mathcal {A}}(x)\Lambda^{-+}(x), \bar \Pi^{-} (x)={\mathcal {A}}(x)\Pi^{-}(x), \bar F^{-+}(x)= {\mathcal {A}}(x)F^{-+}(x)$, $\bar F^{-\,-}(x,y)={\mathcal {A}}(x)F^{-\,-}(y)\tilde{\mathcal {A}}(y), \bar U={\mathcal {A}}(1)U $, one has the $( Z, \bar Y, X)$ system
\begin{align}\label{plant_eq1_new}
Z_t(x,t)&=-\Sigma^+(x)Z_x(x,t) +\Lambda^{++}(x)Z(x,t)+\bar\Lambda^{+-}(x)\bar Y(x,t)+\Pi^+(x) X(t)+\int^x_0 \left[F^{++}(y) Z(y,t)+\bar F^{+-}(y)\bar Y(y,t)\right]dy\\
	\bar Y_t(x,t)&=\Sigma^- (x)\bar Y_x(x,t)+\bar\Lambda^{-\,-} (x)\bar Y(x,t)+\bar\Lambda^{-+}(x) Z(x,t)+\bar\Pi^{-}(x) X(t)+\int^x_0 [\bar F^{-+}(x) Z(y,t)+\bar F^{-\,-}(x, y)\bar Y(y,t)]dy\\
	\dot{X}(t)&=AX(t)+B\bar Y(0, t)
\end{align}
with boundary conditions
\begin{align}
 \bar Y(1,t)=\bar U(t),~Z(0,t)=C\bar Y(0, t)+DX(t).\label{BDm_e_new}
\end{align}}

\subsection{Boundary control law and stability result}\label{main}
For system (\ref{plant_eq1_Hyperbolic})--(\ref{BDm_e_hyperbolic}), the following control law is obtained {according to backstepping transformation~(\ref{tr_eq53}) and the fact $U(t)=\tilde{\mathcal A}(1)\bar Y(1,t)$.}
\begin{eqnarray}\label{eqn-controlaw_new}
	U&\hspace{-2pt}(t)=\hspace{-2pt}& \int_0^1 \tilde{\mathcal A}(1){{K}\left( {1,y} \right) {\mathcal A}(y) Y\left( {y,t} \right)}dy+\int_0^1 \tilde{\mathcal A}(1){{L}\left( {1,y} \right) Z\left( {y,t} \right)}dy+\tilde{\mathcal A}(1)\Phi(1)X(t),\quad \,\,\,
\end{eqnarray}
whose gain kernels are the particular values of the matrices
{
\begin{align}
	K(x, y)\in \mathbb{R}^{n\times n},\Phi(x)\in \mathbb{R}^{n\times \nu}, L(x, y)\in \mathbb{R}^{n\times m}. 
\end{align}}
evaluated at $x=1$. These matrices will be defined in Section~\ref{transf}. Finding $\Phi(x)$ in particular requires setting boundary condition $\Phi(0)$. Define
\begin{align}\label{Phi_delta}
    E_1=A+B\Phi(0).
\end{align}
We can obtain rapid stabilization of the system by choosing $\Phi(0)$ to adequately set the eigenvalues of the $E_1$ matrix, which is always possible due to Assumption~\ref{assump1}~\cite{Wonham1967}, thus obtaining the following result.


\begin{theorem}\label{thm1_new}
	Consider system  (\ref{plant_eq1_Hyperbolic})--(\ref{BDm_e_hyperbolic}), with initial conditions $Z_{0},Y_{0} \in L^2(0,1)$, $X_{0} \in L^2$
	under the  control law (\ref{eqn-controlaw_new}). For all  $C_2>0$ there exists gains $K(1,y)$, $L(1,y)$ and $\Phi(1)$ such that 
 (\ref{plant_eq1_Hyperbolic})--(\ref{BDm_e_hyperbolic}) has a solution $Y(\cdot, t)$, $Z(\cdot,t)\in L^2(0,1)$, $X(t) \in L^2$ for $t>0$, and the following inequality is verified for some $C_1>0$:
	\begin{eqnarray}
		\|Z(\cdot, t)\|^2_{L^2}+\|Y(\cdot, t)\|^2_{L^2}+\|X(t)\|^2_{L^2}
\leq C_1 \mathrm{e}^{-C_2 t}\Big(\|Z_{0}\|^2_{L^2}+\|Y_{0}\|^2_{L^2}+\|X_0\|^2_{L^2}
		\Big).
	\end{eqnarray}
\end{theorem}
The proof of Theorem~\ref{thm1_new} is given in Section \ref{stability}.

\section{Controller Analysis}\label{contdesign}
In this section, we outline the process leading to (\ref{eqn-controlaw_new}). Employing the backstepping method, we begin by introducing the target system in Section~\ref{sec-target}; follow by, the backstepping transformation (of Volterra type) is presented in Section~\ref{transf}; the well-posedness of the kernel equations is established in Theorem~\ref{thm2}.
\subsection{Target system}\label{sec-target}
Inspired by~\cite{Auriol2016}, we design a target system as follows
\begin{align}\label{target10}
	\sigma_t&=\Sigma^-(x)\sigma_x +\Omega(x)\sigma,
	\\
	Z_t&=-\Sigma^+(x)Z_x+\Lambda^{++}(x)  Z+  \Lambda^{+-}(x)\tilde{\mathcal A}(x)\sigma 
	+\Xi_1(x)X+\int^x_0 \Xi_2(x, y)\sigma(y, t)dy
	+\int^x_0 \Xi_3(x, y)Z(y, t)dy\label{target_phi31}, \,\,\,\,\,\\
	\dot{X}&=E_1X+B\sigma(0, t),\label{target-X}
\end{align}
with boundary conditions
\begin{align}
	\sigma(1,t)=0,\,
	Z(0,t)=E_2X(t)+C\sigma(0, t),\label{BDtarg}
\end{align}
where
{
\begin{align}
	\sigma&=\left[
	\begin{array}{c}
		\sigma_1\\
		\sigma_2\\
		\vdots\\
		\sigma_{n}\\
	\end{array}
	\right],
	~E_1=A+B\Phi(0),~E_2={C\Phi(0)+D}.
\end{align}}
and where the values of $\Xi_1(x, y)$, $\Xi_2(x, y)$, and $\Xi_3(x, y)$ are obtained in terms of the inverse backstepping transformation. {The matrix $\Omega(x)$ is a strictly lower triangular matrix }which allows us to solve the matrices of kernel functions line by line, avoiding the multiplications between two unknown kernel functions, given in Section~\ref{transf}. The stability of this target system is shown in Section~\ref{stability}.
\subsection{Backstepping transformation}\label{transf}
To begin with, drawing from \cite{c2} for inspiration, we present a backstepping transformation of Volterra type for $\bar Y$.
\begin{align}
	\sigma =\bar Y-\int_0^x {{K}\left( {x,y} \right)\bar Y\left( {y,t} \right)}dy-\int_0^x {{L}\left( {x,y} \right)Z\left( {y,t} \right)}dy-\Phi(x)X(t)\label{tr_eq53}
\end{align}
where $K(x,y), L(x,y)$ and $\Phi(x)$ verify 
{
\begin{align}
	&\Sigma^-(x) K_{x}(x,y)+K_{y}(x,y)\Sigma^-(y)
	-K(x,y)[\bar\Lambda^{-\,-}(y)-\Sigma^-_y(y)]+L(x,y)\bar\Lambda^{+-}(y)+\Omega(x)K(x,y)+\bar F^{-\,-}(x,y)
    \notag\\
    &-\int^x_y K(x,s)\bar F^{-\,-}(s,y)ds-\int^x_y L(x,s)\bar F^{+-}(y)ds=0,\label{ker_1}\\
	&\Sigma^-(x) L_{x}(x,y)-L_{y}(x,y)\Sigma^+(y)+L(x,y)[\Lambda^{++}(y)-\Sigma^{+}_y(y)]-K(x,y)\bar\Lambda^{-+}(y)+\Omega(x)L(x,y)+\bar F^{-+}(x)\notag\\
    &-\int^x_y K(x,s)\bar F^{-+}(s)ds-\int^x_y L(x,s) F^{++}(y)ds=0,\label{ker_2}\\
	&\Sigma^-(x)\Phi_x(x)-\Phi(x) A+\bar\Pi^-(x)+\Omega(x)\Phi(x)+L(x, 0)\Sigma^+(0) D-\int^x_0 [K(x,y)\bar \Pi^-(y)+L(x,y)\Pi^+(y)]dy=0,\label{ker_phi}\\
    &\Sigma^-(x) L(x,x)+L(x,x)\Sigma^+(x)+\bar\Lambda^{-+}(x)=0,\label{ker_4} \\
	&\Sigma^-(x) K(x,x)-K(x,x)\Sigma^-(x)+\bar\Lambda^{-\,-}(x)-\Omega(x)=0,\label{ker_3}\\
	&K(x, 0)\Sigma^-(x)+L(x, 0)\Sigma^+(0) C-\Phi(x) B=0,\label{ker_5}
 \end{align}}
 with
\begin{equation}
	\Omega(x)=\left[
	\begin{array}{cccc}
		0&0&0&0\\
		\omega_{21}(x)&0&0&0\\
		\vdots&\ddots&0&0\\
		\omega_{n1}(x)&\cdots&\omega_{nn-1}(x)&0\\
	\end{array}
	\right].
\end{equation}
{
Developing (\ref{ker_1})--(\ref{ker_phi}), we have the following set of kernel PDEs:
\begin{align}
    \Sigma^-_{ii}(x)\partial_xK_{ij}(x,y)+\partial_yK_{ij}(x,y)\Sigma^-_{jj}(y)=&\Sigma^n_{k=1}\{K_{ik}(x,y)[\bar\Lambda^{-\,-}_{kj}(y)-\partial_y\Sigma^-_{kj}(y)]-L_{ik}(x,y)\bar\Lambda^{+-}_{kj}(y)-\Omega_{ik}(x)K_{kj}(x,y)\}\notag\\
    &+\Sigma^n_{k=1}\{-\bar F^{-\,-}_{ij}(x,y)+\int^x_y K_{ik}(x,s)\bar F^{-\,-}_{kj}(s,y)ds+\int^x_y L_{ik}(x,s)\bar F^{+-}_{kj}(y)ds\},\label{kernel_K_Dev}\\
    \Sigma^-_{ii}(x) \partial_x L_{i\tau}(x,y)-\partial_y L_{i\tau}(x,y)\Sigma^+_{\tau\tau}(y)=&\Sigma^{m}_{k=1}\{-L_{i k}(x,y)[\Lambda^{++}_{k \tau}(y)-\partial_y\Sigma^{+}_{k \tau}(y)]+K_{i k}(x,y)\bar\Lambda^{-+}_{k\tau}(y)-\Omega_{i k}(x)L_{k\tau}(x,y)\}\notag\\
    &+\Sigma^{m}_{k=1}\{-\bar F^{-+}_{i \tau}(x)+\int^x_y K_{i k}(x,s)\bar F^{-+}_{k\tau}(s)ds+\int^x_y L_{i k}(x,s) F^{++}_{k\tau}(y)ds \},\label{kernel_L_Dev}\\
    \Sigma^-_{ii}(x)\partial_x\Phi_{i\xi}(x)=&\Sigma^{\nu}_{k=1}\Phi_{ik}(x) A_{k\xi}-\bar\Pi^-_{i \xi}(x)-\Omega_{ik}(x)\Phi_{k\xi}(x)-L_{ik}(x, 0)[\Sigma^+(0) D]_{k\xi}\notag\\
    &+\Sigma^{\nu}_{k=1}\int^x_0 [K_{ik}(x,y)\bar \Pi^-_{k\xi}(y)+L_{ik}(x,y)\Pi^+_{k\xi}(y)]dy,\label{kernel_Phi_Dev}
\end{align}
with the set of kernel boundary conditions:
\begin{align}
L_{i\tau}(x,x)=&\frac{\bar\Lambda^{-+}_{i\tau}(x)}{\Sigma^-_{ii}(x)+\Sigma^+_{\tau\tau}(x)},~~\text{for}~~1\le i\le n, 1\le \tau \le m,\label{ker_4_dev} \\
K_{ij}(x,x)=&\frac{-\bar\Lambda^{-\,-}_{ij}(x)+\Omega_{ij}(x)}{\Sigma^-_{ii}(x)-\Sigma^-_{jj}(x)}, ~~\text{for}~~1\le i,j\le n ~(\Sigma^-_{ii}\neq \Sigma^-_{jj}),\label{ker_3_dev}\\
K_{ij}(x, 0)=&\frac{\Sigma^{n}_{k=1}-L_{ik}(x, 0)[\Sigma^+(0) C]_{kj}+\Phi_{ik}(x) B_{kj}}{\Sigma^-_{jj}(x)}, ~~\text{for}~~1\le i,j\le n.\label{ker_4_dev}
\end{align}}

The kernel equations are derived conventionally through a meticulous yet methodical process involving differentiation in the transformation, substitution of the original and target equations, and integration by parts. However, for conciseness, the intricate details are omitted. {The introduction of $\Omega(x)$ is to guarantee that the piecewise kernels $K_{ij}$ all have one and only one boundary condition in each sectional area and avoid over- or under-determined cases. They are chosen as } $\omega_{i,j}(x)=(\Sigma^-_{ii}(x)-\Sigma^-_{jj}(x))K_{ij}(x, x)+\bar\Lambda^{-\,-}_{ij}(x),~i>j,~i=2,\cdots,n$. Notice that for $i,j$ belonging to the same block, $\Sigma^-_{ii}(x)=\Sigma^-_{jj}(x)$ but also $\bar\Lambda^{-\,-}_{ij}(x)=0$, thus resulting in $\omega_{ij}(x)=0$. Therefore there are no singularities in (\ref{ker_3})  for $i,j$ such that $\Sigma^-_{ii}(x)=\Sigma^-_{jj}(x)$. Depending on the slopes of the characteristics of the kernel equations, one or two boundary conditions may be needed, and the structure of $\Omega(x)$ allows to accommodate the kernels in (\ref{ker_3}) that only accept one boundary condition.

The structure of the kernel equations is similar to \cite{Auriol2016}. For $n\geq i\geq 2$, the kernel equations for $K_{ij}$, $L_{ij}$ and $\Phi_{ij}$ seem to be nonlinear { since $\Omega(x)$ is the function of $K(x,x)$ which results in products of undetermined functions if we try to solve all the kernels synchronously}. However, one can start by solving  $K_{1j}$, $L_{1j}$ and $\Phi_{1j}$, which are linear and can be proven solvable. Then, they become known coefficients of the equations verified $K_{2j}$, $K_{2j}$ and $\Phi_{2j}$. Thus $K_{2j}$, $L_{2j}$ and $\Phi_{2j}$ become also linear and solvable. In the same recursive manner~\cite{Auriol2016}, we can obtain the solution of each kernel equation. Regarding the well-posedness of $K(x,y), L(x,y)$,
the following result holds.

\begin{theorem}\label{thm2}
	There exists a unique bounded solution to the kernel equations (\ref{ker_1})--(\ref{ker_5}), namely
	$K(x,y)$, $L(x, y)$, $\Phi(x)$; in particular, there exists  positive numbers $\mathcal N, \mathcal M$ such that
	\begin{eqnarray}
		\Vert K(x,y) \Vert_{\infty},\Vert L(x,y)\Vert_{\infty} , \Vert \Phi(x)\Vert_{\infty} \leq \mathcal N \mathrm{e}^{\mathcal Mx}.
	\end{eqnarray} 
\end{theorem}

The proof follows the approach outlined in~\cite{c1}, and is omitted here. It relies on utilizing the method of characteristics to reformat (\ref{ker_1})--(\ref{ker_5}) into integral equations, followed by establishing a solution using a series of successive approximations, whose convergence is recursively validated. It's evident that the procedures of~\cite{c1} can be readily adjusted to accommodate the presence of blocks, integral terms, and discrepancies in boundary conditions with minimal difficulty.

Since the kernels in (\ref{tr_eq53}) exhibit bounded behavior, the transformation proves invertible according to the principles of Volterra integral equations. Consequently, one can define
\begin{align}
	\bar Y &=\sigma +\int_0^x {{{\mathord{\buildrel{\lower3pt\hbox{$\scriptscriptstyle\smile$}}
					\over K} }}\left( {x,y} \right)\sigma\left( {y,t} \right)} dy+\int_0^x {{{\mathord{\buildrel{\lower3pt\hbox{$\scriptscriptstyle\smile$}}
					\over L} }}\left( {x,y} \right)Z\left( {y,t} \right)}dy+{{\mathord{\buildrel{\lower3pt\hbox{$\scriptscriptstyle\smile$}}
				\over \Phi} }}\left( x \right)X.\label{eq9}
\end{align}
with bounded kernels. The transformatio and its inverse map {$(L^2)^{n+m}\times \mathbb R^{\nu}$ functions into $(L^2)^{n+m}\times \mathbb R^{\nu}$ functions and back}. (see e.g.~\cite{c6}).

Derived from the inverse transformation, the kernels $\Xi_1(x)$, $\Xi_2(x, y)$, $\Xi_3(x, y)$ appearing in (\ref{target_phi31}) are
\begin{align}
	\Xi_1(x)&=\Lambda^{+-}(x)\tilde{\mathcal A}(x){{{\mathord{\buildrel{\lower3pt\hbox{$\scriptscriptstyle\smile$}} \over \Phi} }}}(x)+\Pi^+(x)\int^x_0 F^{+-}(y)\tilde{\mathcal A}(y){{{\mathord{\buildrel{\lower3pt\hbox{$\scriptscriptstyle\smile$}} \over \Phi} }}}(y)dy,\\
	\Xi_2(x, y)&=\Lambda^{+-}(x)\tilde{\mathcal A}(x){{\mathord{\buildrel{\lower3pt\hbox{$\scriptscriptstyle\smile$}} \over K} }}(x, y)+F^{+-}(x)\tilde{\mathcal A}(y)+\int^x_0 F^{+-}(y)\tilde{\mathcal A}(y) {{\mathord{\buildrel{\lower3pt\hbox{$\scriptscriptstyle\smile$}} \over K} }}(s, y)ds,\\
	\Xi_3(x, y)&=\Lambda^{+-}(x)\tilde{\mathcal A}(x) {{\mathord{\buildrel{\lower3pt\hbox{$\scriptscriptstyle\smile$}} \over L} }}(x, y)+F^{++}(x)+\int^x_0 F^{+-}(y) \tilde{\mathcal A}(y){{\mathord{\buildrel{\lower3pt\hbox{$\scriptscriptstyle\smile$}} \over L} }}(s, y)ds.
\end{align}
from which it can be deduced  that they are bounded kernels by the method of successive approximations.

\section{Stability and Analysis of Closed Loop}\label{stability}

This section proves Theorem~\ref{thm1_new}. First, in Section~\ref{sec-solution}, the solution of  (\ref{target10})--(\ref{BDtarg}) is studied with the method of characteristics. This helps to find stability conditions in Section~\ref{sec-solution}. Then, a Lyapunov analysis in Section~\ref{sec-lyap} shows exponential stability.
\subsection{A semi-explicit solution for the target system}\label{sec-solution}
We begin by solving (\ref{target10})--(\ref{BDtarg}) with the method of characteristics. Writing down the solution {for $\sigma(x, t)=\left[\sigma_1~~\sigma_2~~\cdots~~\sigma_{n}\right]^T $:
	\begin{align}
		\sigma_i&=\left\{\begin{aligned} &\sigma_i(x+\Sigma^{-}_{ii}(x)t,0)+ \sum_{j=1}^{i-1}\int_0^t \omega_{ij}\left(x+\Sigma^-_{ii}(x)(t-s)\right)\sigma_j\left(x+\Sigma^-_{ii}(t-s),s \right) ds, 0\leq t \leq t_i(x),\\
&0,~~~~~~~~~~~~~~~~~~~~~~~~~~~~~~t>t_i(x),~i=1,2,\cdots,n
		\end{aligned} \right.
	\end{align}
	where
	\begin{align}
		t_i(x)=\frac{(1-x)}{\Sigma^-_{ii}(x)},~i=1,2,\cdots,n.
	\end{align}}
	Thus, $\sigma(x,t)$ converges to zero in finite time $\frac{1}{\Sigma^-_{nn}(x)}$.
For $t>\frac{1}{\Sigma^-_{nn}(x)}$,
\begin{align}\label{eq:explicit}
	Z_t(x, t)=&-\Sigma^+(x) Z_x(x, t)+\Lambda^{++}(x) Z(x, t)+\Xi_1(x)X{+}\int^x_0 \Xi_3(x, y)Z(y, t)dy,\\
	\dot{X}=&E_1X,\label{dotxeqn}\\
	Z(0,t)=&BX.
\end{align}
Solving for $X$ we get $X(t)=\mathrm{e}^{E_1 t}X(0)$, where we have employed the matrix exponential. Then
\begin{align}\label{eq:semi-explicit}
	Z_t(x, t)=&-\Sigma^+ Z_x(x, t)+\Lambda^{++} Z(x, t)+\Xi_1(x)X(0)\mathrm{e}^{E_1 t}+\int^x_0 \Xi_3(x, y)Z(y, t)dy,\\
	Z(0,t)=&B\mathrm{e}^{E_1 t}X(0).
\end{align}
Using the method of characteristics, Volterra-type integral equations can be derived for the components of $Z$. While the specifics are omitted, a unique $L^2$ solution for $Z$ can always be obtained which ensures at least stability of the target system.

Stability hinges solely on the Hurwitz property of $E_1$ for ensuring the exponential stability of the state's origin of (\ref{eq:explicit}). However, for rapid arbitrary stabilization, it's essential to configure the eigenvalues of $E_1=A+B\Phi(0)$, possibly through pole placement. Consequently, if we opt for boundary conditions $\Phi(0)$ such that
\begin{align}\label{E_1condition}
      E_1+E_1^T<-2c I.
\end{align}
then the zero equilibrium of $X$ in (\ref{dotxeqn}) is exponentially stable with a convergence rate of at least {$c>1$}. This is always possible to achieve by Assumption~\ref{assump1}.

\subsection{Lyapunov-based stability analysis of target system}\label{sec-lyap}
Subsequently, we employ a Lyapunov function to analyze the stability of the target system, demonstrating exponential stability of the origin with a predetermined convergence rate. Define
\begin{align}
	V=&{\zeta_1} X^T X+{\zeta_2}\int^1_0 \mathrm{e}^{\delta x}\sigma^T(x, t) (\Sigma^-(x))^{-1}\sigma(x, t)dx+ \int^1_0 \mathrm{e}^{-\delta x}Z^T(x, t) (\Sigma^+(x))^{-1}Z(x, t)dx,\label{eq54}
\end{align}
{where the parameters $\zeta_1, \zeta_2$, $\delta$ belong to $\mathbb{R}$.} Taking the derivative of {(\ref{eq54})} with respect to $t$, we obtain
\begin{align}\label{DLy}
	\dot V=&2{\zeta_1} X^T \dot X+2{\zeta_2}\int^1_0 \mathrm{e}^{\delta x}\sigma^T(x, t) (\Sigma^-(x))^{-1}\sigma_t(x, t)dx+ 2\int^1_0 \mathrm{e}^{-\delta x}Z^T(x, t) (\Sigma^+(x))^{-1}Z_t(x, t)dx
\end{align}
By substituting equations (\ref{target10})–(\ref{BDtarg}) into equation (\ref{DLy}) and subsequently employing integration by parts, we arrive at
\begin{eqnarray}
\dot{V} &=&-2\zeta_1c X^TX+2\zeta_1X^T B\sigma(0, t)-\zeta_2\sigma^T(0, t)\sigma(0, t)-\zeta_2 \int^1_0 \mathrm{e}^{\delta x}\sigma^T(x, t) (\delta I-2(\Sigma^{-}(x))^{-1}\Omega(x)) \sigma(x, t)dx\nonumber\\
&&- \int^1_0 \mathrm{e}^{-\delta x}Z^T(x, t)(\delta I+2(\Sigma^+(x))^{-1}\Lambda^{++}(x))  Z(x, t)dx-2\int^1_0 \mathrm{e}^{-\delta x} Z^T(x, t)(\Sigma^+(x))^{-1}\Lambda^{+-}(x) \tilde{\mathcal A}(x)\sigma(x,t)dx\nonumber\\
&&-2\int^1_0 \mathrm{e}^{-\delta x} Z^T(x, t)(\Sigma^+(x))^{-1}\Xi_1 (x)Xdx-2\int^1_0 \mathrm{e}^{-\delta x} Z^T(x, t)(\Sigma^+(x))^{-1}\int^x_0 \Xi_2 (x, y)\sigma(y, t)dydx\nonumber\\
&&-2\int^1_0 \mathrm{e}^{-\delta x} Z^T(x, t)(\Sigma^+(x))^{-1}\int^x_0 \Xi_3 (x, y)Z(y, t)dydx+ Z^T(0,t)Z(0, t)
. \,\,\,\,\,\, \label{eq60}
\end{eqnarray}
Concerning the final line of equation (\ref{eq60}), {utilizing $Z(0,t)=E_2X+C\sigma(0, t)$, we obtain $Z^T(0,t)Z(0, t)= X^TE^T_2E_2X+2 X^TE^T_2C\sigma(0, t)+\sigma^T(0, t)C^TC\sigma(0, t)$}. {Using Young's inequality}, the initial and final lines of (\ref{eq60}) {are transformed to}:{
\begin{eqnarray}
&&-X^T(2\zeta_1c I- E^T_2E_2)X+2X^T(\zeta_1 B+ E^T_2C) \sigma(0, t)-\sigma^T(0, t)(\zeta_2- C^TC)\sigma(0, t) \nonumber\\
&\leq & -(2\zeta_1c-M_1) X^TX- (\zeta_2-M_2) \sigma^T(0, t)\sigma(0, t), \label{b66-2}
 \end{eqnarray}
  where $M_1=\Vert E_2 \Vert^2+1$, $M_2=\Vert \zeta_1 B+ E^T_2C \Vert^2+\Vert C \Vert^2$.}
{Applying Young's inequality for twice, we have:}
\begin{eqnarray}
&&-2\int^1_0 \mathrm{e}^{-\delta x} Z^T(x, t)(\Sigma^+)^{-1}(x)\Lambda^{+-}(x) \tilde{\mathcal A}(x)\sigma(x,t)dx-2\int^1_0 \mathrm{e}^{-\delta x} Z^T(x, t)(\Sigma^+)^{-1}(x)\Xi_1 (x)Xdx\nonumber\\
&\leq &
(1+M_4)\int^1_0 \mathrm{e}^{-\delta x} Z^T(x, t) Z(x,t) dx+M_3 \int^1_0 \mathrm{e}^{\delta x} \sigma^T(x, t) \sigma(x,t) dx
+  X^TX,\label{b67} 
 \end{eqnarray}
where $M_3=\max_{x\in[0,1]}\Vert (\Sigma^+(x))^{-1}\Lambda^{+-}(x)\tilde{\mathcal A}(x) \Vert ^2 $ and $M_4=\max_{x\in[0,1]}\Vert (\Sigma^+(x))^{-1} \Xi_1 (x) \Vert^2$.
 {The last term in the third line of (\ref{eq60}) can be bounded using Young's inequality again as follows:}
\begin{eqnarray}
&&2\int^1_0 \mathrm{e}^{-\delta x} Z^T(x, t)(\Sigma^+(x))^{-1}\int^x_0 \Xi_2(x, y)\sigma(y, t)dydx\nonumber\\
&\leq &
2\int^1_0 \int^1_0 \mathrm{e}^{-\delta x} \vert Z^T(x, t)\vert (\Sigma^+(x))^{-1}  \vert \Xi_2 (x, y) 
\vert\vert\sigma(y, t)\vert dydx\nonumber\\
&\leq &
 \int^1_0 \mathrm{e}^{-\delta x} Z^T(x, t)   Z(x, t)dx+M_5 \int^1_0 \mathrm{e}^{\delta x} \sigma^T(x, t) \sigma(x,t) dx,
 \end{eqnarray}
  where $M_5= \max_{x,y\in[0,1]}\Vert({\Sigma^+(x))^{-1} }\Xi_2 (x, y)  \Vert^2$.
 Finally, {the integral term in the fourth lineof (\ref{eq60}) is also bounded following Cauchy-Schwarz inequality} and satisfies
\begin{eqnarray}
&&2\int^1_0 \mathrm{e}^{-\delta x} Z^T(x, t)(\Sigma^+(x))^{-1}\int^x_0 \Xi_3 (x, y)Z(y, t)dydx\nonumber\\
&\leq &
2\int^1_0\mathrm{e}^{-\delta x/2} \vert Z^T(x, t)\vert (\Sigma^-(x))^{-1}\int^x_0  \mathrm{e}^{-\delta x/2}  \vert \Xi_3 (x, y)\vert\vert Z(y, t)\vert dydx\nonumber\\
&\leq &
2\int^1_0 \int^1_0 \mathrm{e}^{-\delta x/2} \vert Z^T(x, t)\vert (\Sigma^+(x))^{-1}
\mathrm{e}^{-\delta y/2}  \vert \Xi_3 (x, y) 
\vert\vert Z(y, t)\vert dydx\nonumber\\
&\leq & M_6
 \int^1_0 \mathrm{e}^{-\delta x} Z^T(x, t) Z(x, t)dx,
 \label{b69}
 \end{eqnarray}
 with $M_6= 1+\max_{x,y\in[0,1]}\Vert {(\Sigma^+(x))^{-1}}\Xi_3 (x, y) \Vert^2$. Therefore
 \begin{eqnarray}
\dot V&\leq&
 -(2\zeta_1c-M_1-1) X^TX- (\zeta_2-M_2) \sigma^T(0, t)\sigma(0, t) \nonumber
 \\ &&
- \zeta_2\int^1_0 \mathrm{e}^{\delta x}\sigma^T(x, t) ( \delta I-2 (\Sigma^-(x))^{-1}\Omega(x)-M_3I-M_5I)\sigma(x, t)dx
 \nonumber\\ 
&&-\int^1_0 \mathrm{e}^{-\delta x}Z^T(x, t)(\delta I+2(\Sigma^{+}(x))^{-1}\Lambda^{++}-2I-M_4I-M_6I)Z(x, t)dx.\label{eq70}
\end{eqnarray}
 By selecting {$c=(c'+1)/2$ where $c'>0$, $\zeta_1=M_1+1$, $\delta >\max \{(2 {\Vert \Lambda^{++}(x) \Vert+c')\Vert(\Sigma^+(x))^{-1}}\Vert +M_6+2+M_4,(c'+2 \max_{x\in[0,1]} {\Vert \Omega(x) \Vert)\Vert(\Sigma^-(x))^{-1}\Vert} +1\}$, and $\zeta_2 >\max\{M_3+M_5,M_2\}$}, we obtain
\begin{align}
	\dot V\leq&
	-c'\zeta_1 X^TX
	-c' \zeta_2 \int^1_0 \mathrm{e}^{\delta x}\sigma^T(x, t)  (\Sigma^-(x))^{-1}\sigma(x, t)dx\nonumber
	\\
	&-c' \int^1_0 \mathrm{e}^{-\delta x}Z^T(x, t)(\Sigma^+)^{-1}  Z(x, t)dx \leq -c' V
\end{align}
Thus adequately setting $\Phi(0)$ to verify (\ref{E_1condition}), 
an arbitrary convergence rate $c'>0$ is achieved for $V$.

\section{Rapid Stabilization of {\it N}-layer Timoshenko Beam}\label{sec-extension}

Here, we apply the method developed in Section 3 through Section 5 to the {\it N}-layer Timoshenko beam for achieving rapid stabilization with arbitrarily fast decay rate. {From the matrices of (\ref{Tm_e1})-(\ref{Tm_e2}), we notice that $\Pi= \Pi^{-\,-} = \Pi^{-+}=-\Pi^{++}=-\Pi^{+-}$ which is much simpler than the general case in~(\ref{plant_eq1_Hyperbolic}) and (\ref{plant_eq2_Hyperbolic})}. The main results are as follows:


Applying the procedure in Section \ref{sec2} through Section \ref{stability} to the system (\ref{Tm_e1})-(\ref{TBDm_e}), we obtain the following control law
\begin{align}\label{eqn-controlaw}
	U^{pq}&= \int_0^1 \tilde{\mathcal A}(1){{K}\left( {1,y} \right) {\mathcal A}(y) Y\left( {y,t} \right)}dy+\int_0^1 \tilde{\mathcal A}(1){{L}\left( {1,y} \right) Z\left( {y,t} \right)}dy+\tilde{\mathcal A}(1)\Phi(1)X,
\end{align}
whose gain kernels are found from (\ref{ker_1})-(\ref{ker_5}).
Now, the control laws of the {$N$}-layer Timoshenko beams can be expressed as follows:
\begin{align}
	U^{or}_{i}(t)&=\left\{\begin{array}{ll} \frac{U^{pq}_i}{\sqrt{\eta^{or}_i}} -\Sigma_{ii}v^{or}_{i,t}(1,t), & {i~\mathrm{is~odd}},\\
			\frac{U^{pq}_i}{\sqrt{\alpha^{or}_i}}-\Sigma_{ii}\theta^{or}_{i,t}(1,t),&  {i~\mathrm{is~even}}
		\end{array} \right.\notag\\
  i&=1,2,\cdots,2N.\label{U_1}
\end{align}
where $U^{or}_{i}(t)$ are the ordered control inputs. $v^{or}_{i,t}(1,t)$ and $\theta^{or}_{i,t}(1,t)$ are the ordered boundary states.

{Theorem~\ref{thm1_new} can be applied to the multilayer Timoshenko beam plant with the main result stated in terms of the beam states.}
\begin{theorem}\label{thm1}
	Consider system (\ref{plant_eq1})--(\ref{plant_bc2}), with initial conditions $v_{i0},\theta_{i0} \in H^1(0,1), v_{i0t},\theta_{i0t}\in L^2,~i=1,2,\cdots,N$,
	under the  control law (\ref{U_1}). 
 For all $C_2>0$ there exists gains $K(1,y)$, $L(1,y)$ and $\Phi(1)$ such that 
 (\ref{plant_eq1})--(\ref{plant_bc2}) has a solution  $v_i(\cdot, t)$, $\theta_i(\cdot, t)\in H^1(0,1)$, $u_{it}(\cdot, t)$, $\theta_{it}(\cdot, t)\in L^2(0,1)$ for $t>0$, and the following inequality is verified for some $C_1>0$
	\begin{eqnarray}
		&&\sum_{i=1}^{N}\|v_i(\cdot, t)\|^2_{H^1}+\|\theta_i(\cdot, t)\|^2_{H^1}+\|v_{it}(\cdot, t)\|^2_{L^2}+\|\theta_{it}(\cdot, t)\|^2_{L^2}\nonumber\\
		\nonumber\\
		&\leq& C_1 \mathrm{e}^{-C_2 t}\Big(\sum_{i=1}^{N}\|v_{i0}\|^2_{H^1}+\|\theta_{i0}\|^2_{H^1}+\|v_{i0t}\|^2_{L^2}+\|\theta_{i0t}\|^2_{L^2}
		\Big).
	\end{eqnarray}
\end{theorem}

The results follows by combining Theorem~\ref{thm1_new} with the transformation (\ref{eq2})--(\ref{eq2-bis}) and setting $\Phi(0)$ so that  $E_1=A+B\Phi(0)$ verifies $C_2=(c'+1)/2-1=c-1$.
\begin{figure*}[t]
\begin{centering}
		\includegraphics[width=15cm]{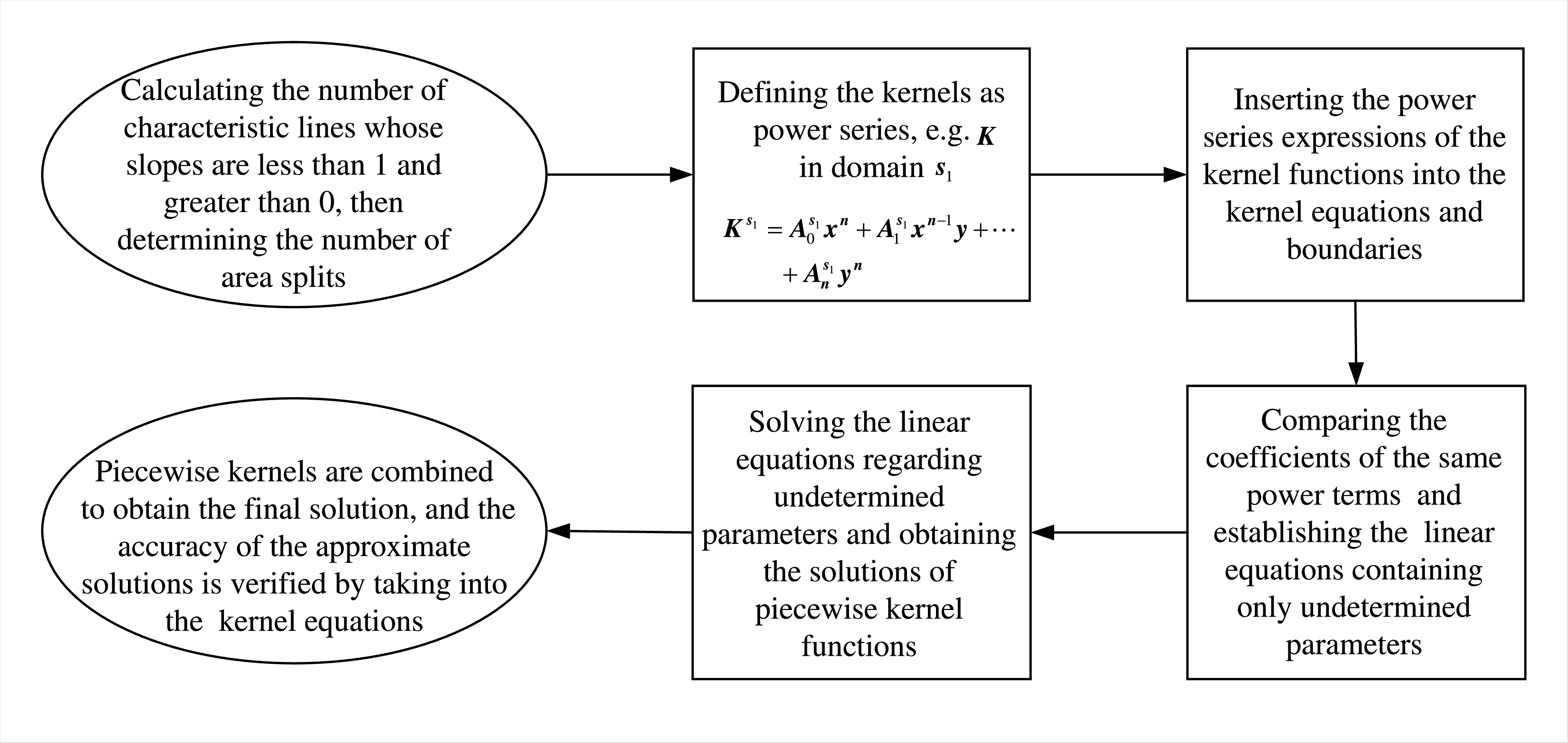}
		\caption{The procedure of power series method}\label{Fig-power}
\end{centering}
\end{figure*} 
\begin{figure*}
\begin{centering}
		\includegraphics[width=8cm]{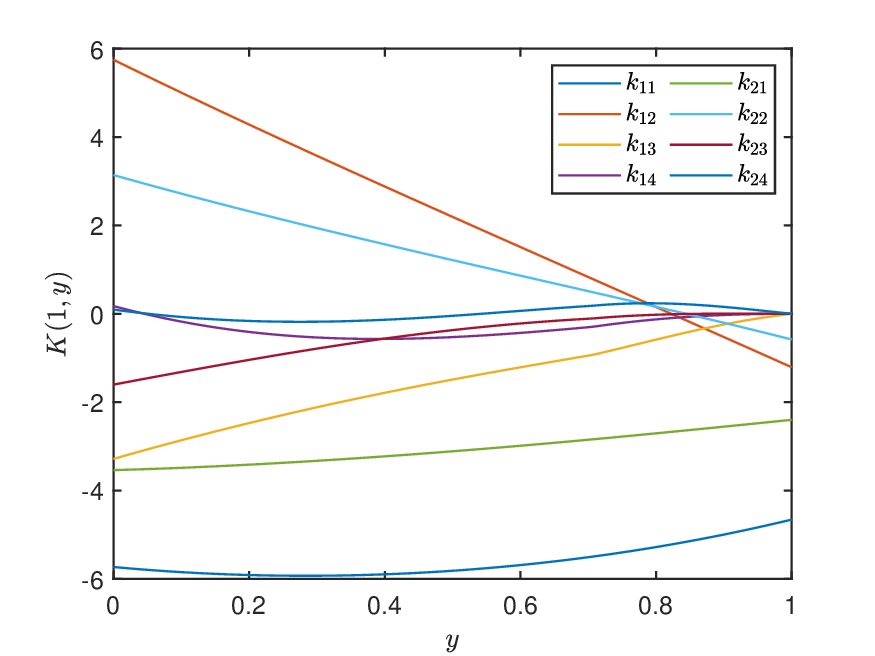}
		\includegraphics[width=8cm]{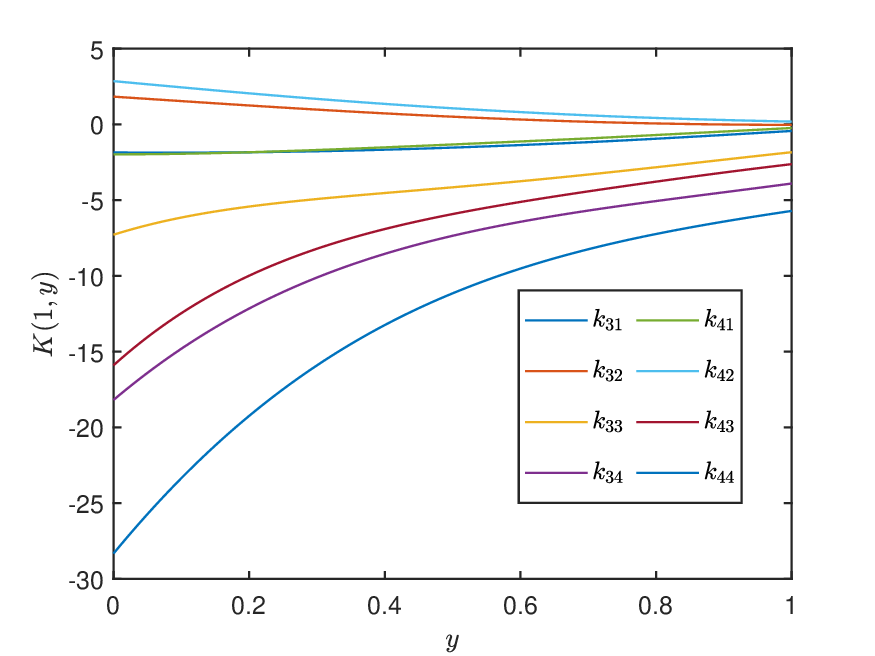}\\
		\includegraphics[width=8cm]{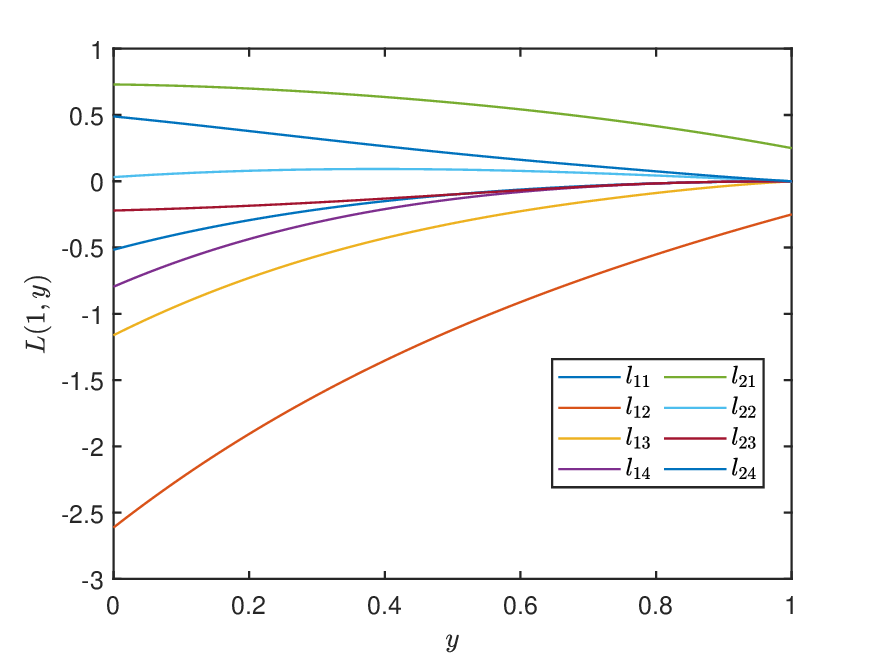}
		\includegraphics[width=8cm]{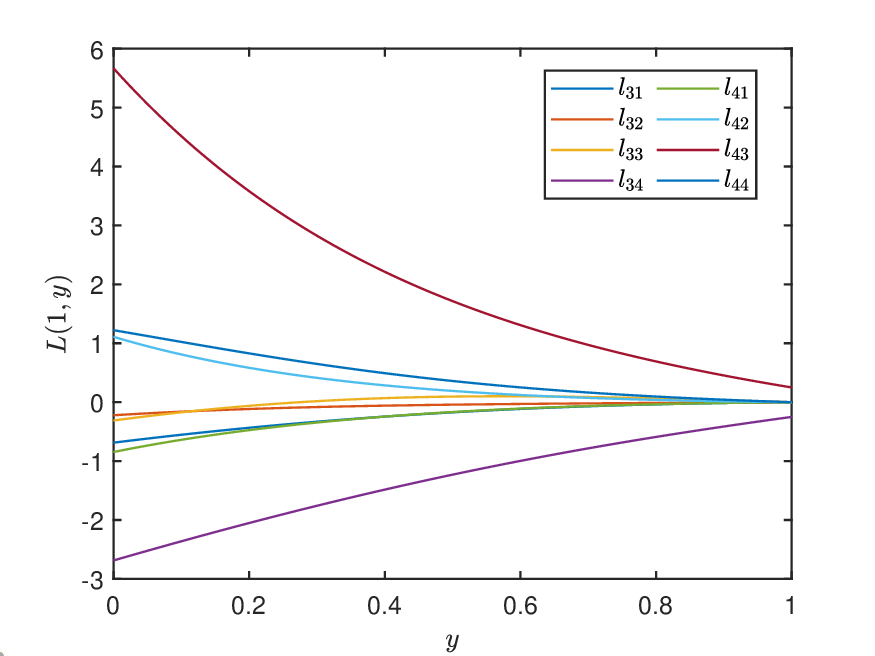}
		\caption{Solutions of Timoshenko gain kernels $K_{ij}(1,y),L_{ij}(1,y), 1\le i\le4,1\le j\le4$ (from left to right).}
		\label{Fig-kernels}
\end{centering}
\end{figure*} 

\begin{figure*}
	\begin{centering}
		\includegraphics[width=8cm]{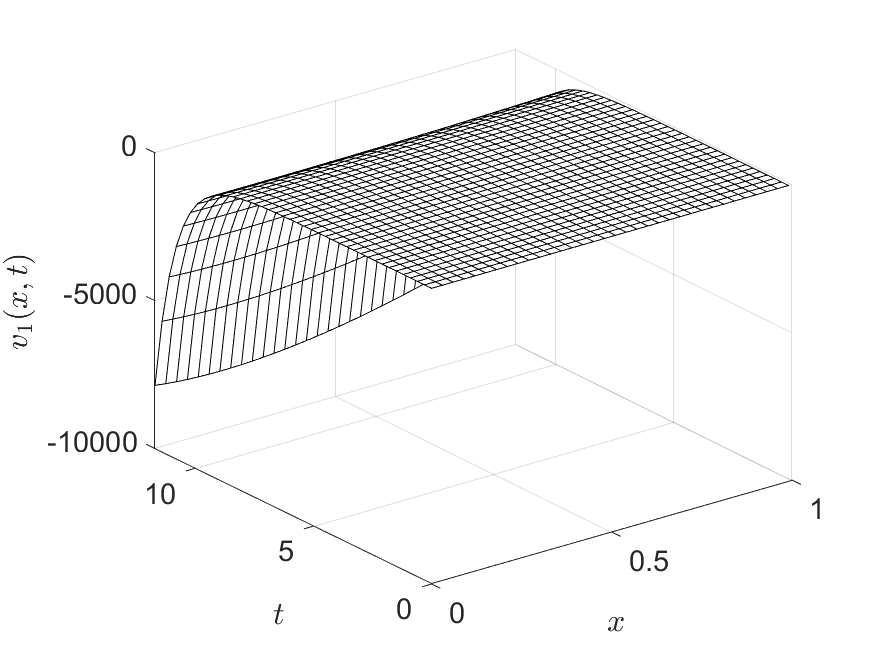}
		\includegraphics[width=8cm]{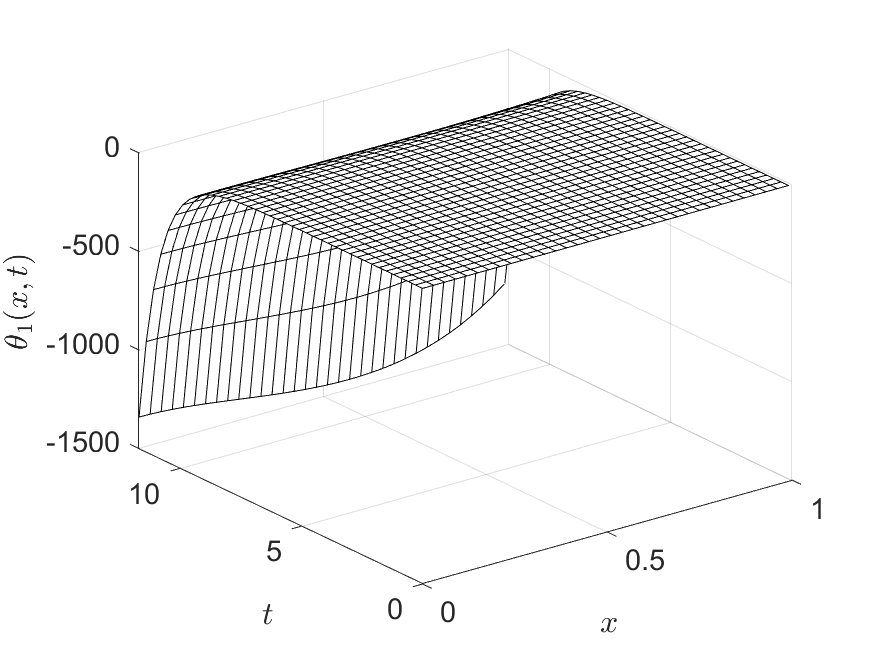}\\
		\includegraphics[width=8cm]{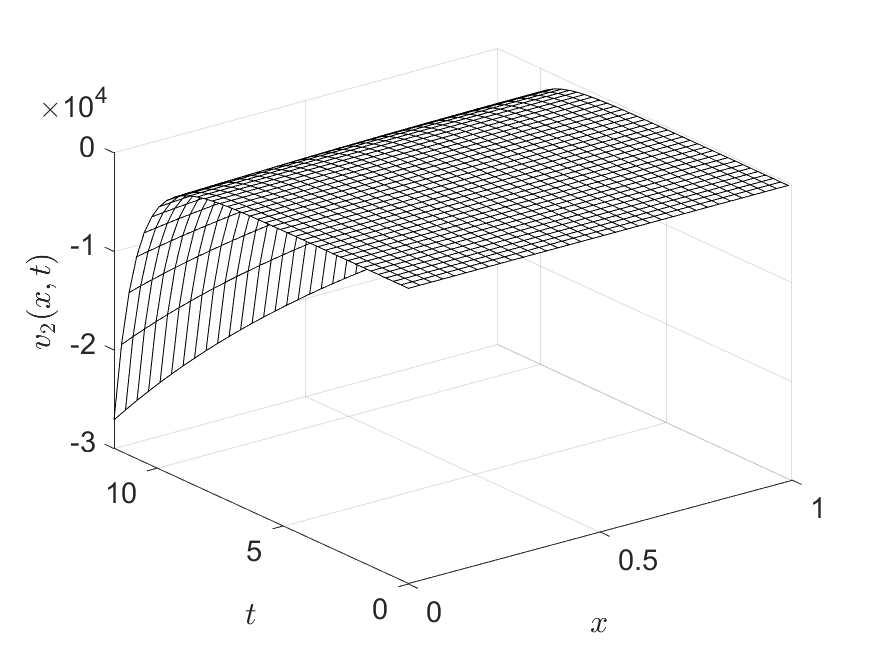}
		\includegraphics[width=8cm]{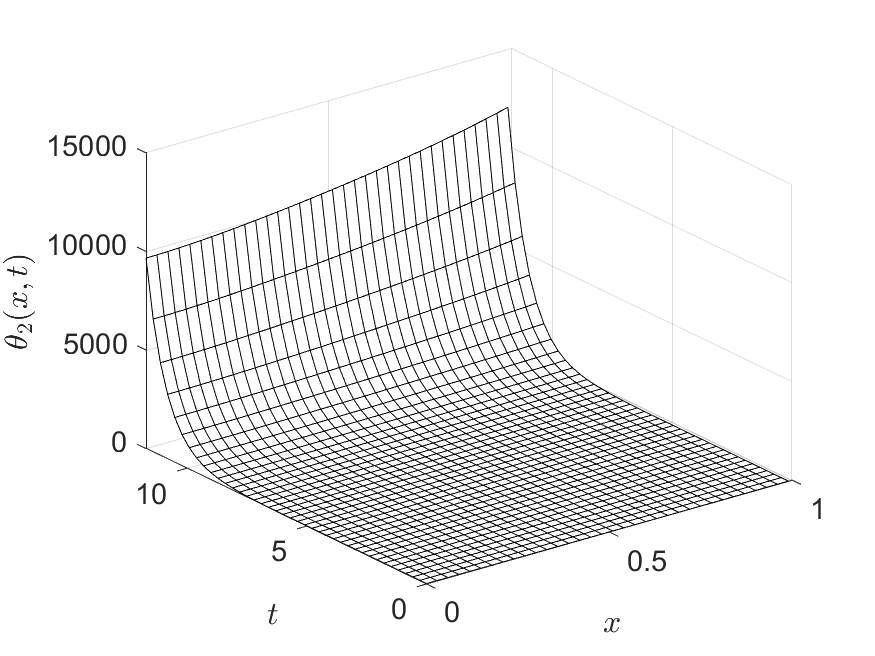}
		\caption{Evolution of open-loop Timoshenko states $v_1(x,t),\theta_1(x,t), v_2(x,t),\theta_2(x,t)$ (from left to right).}
		\label{Fig-open}
	\end{centering}
\end{figure*} 

\begin{figure*}
	\begin{centering}
		\includegraphics[width=8cm]{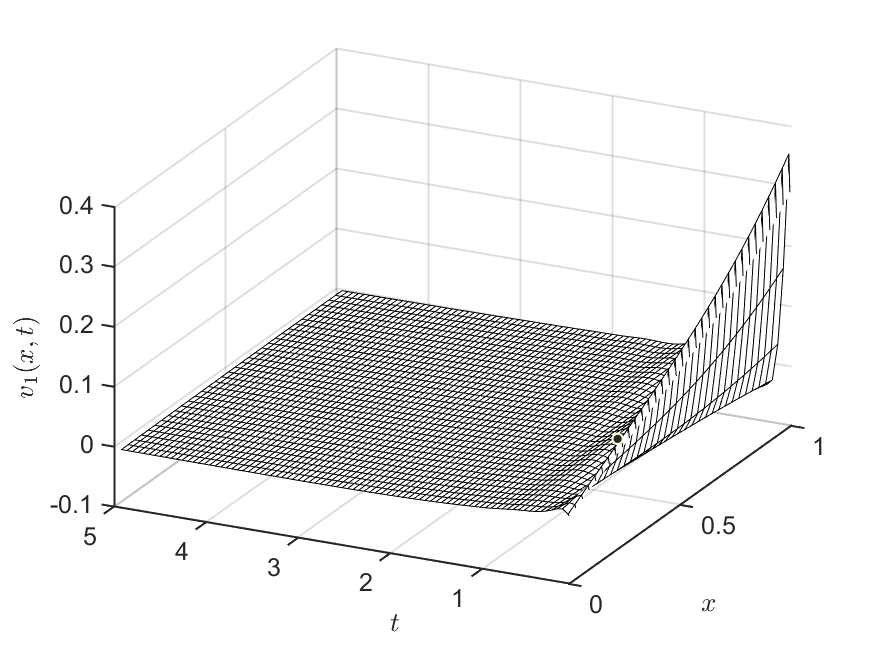}
		\includegraphics[width=8cm]{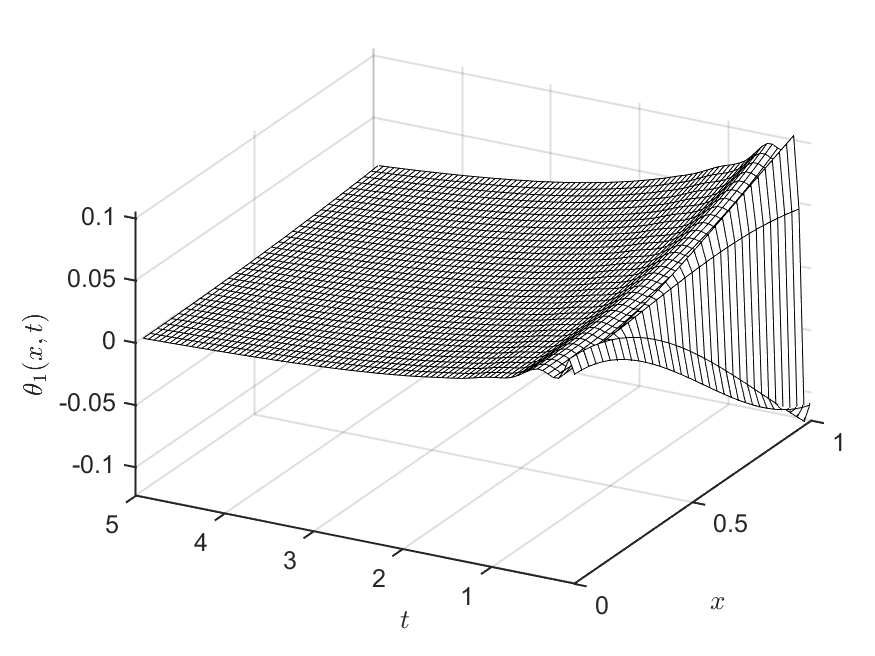}\\
		\includegraphics[width=8cm]{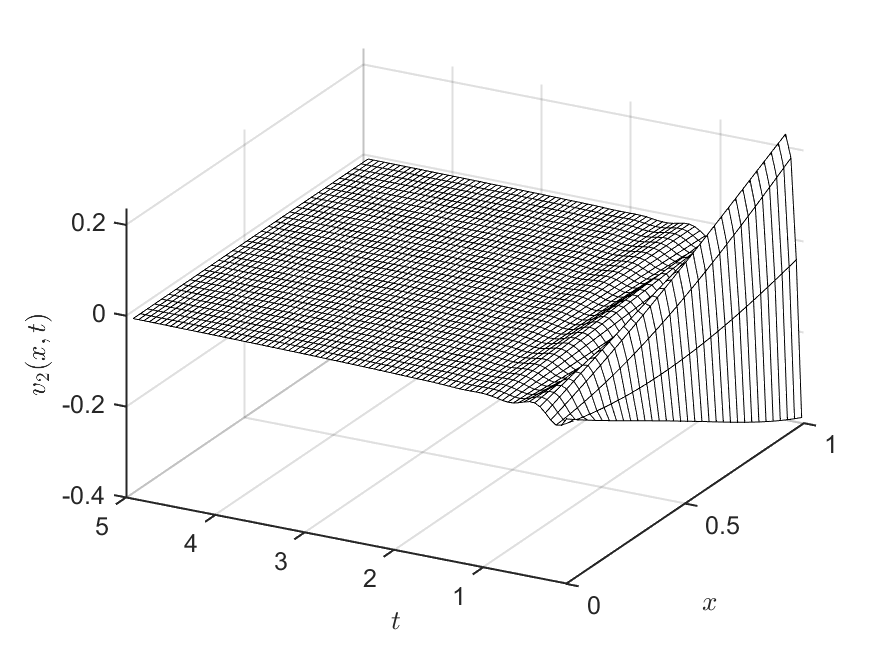}
		\includegraphics[width=8cm]{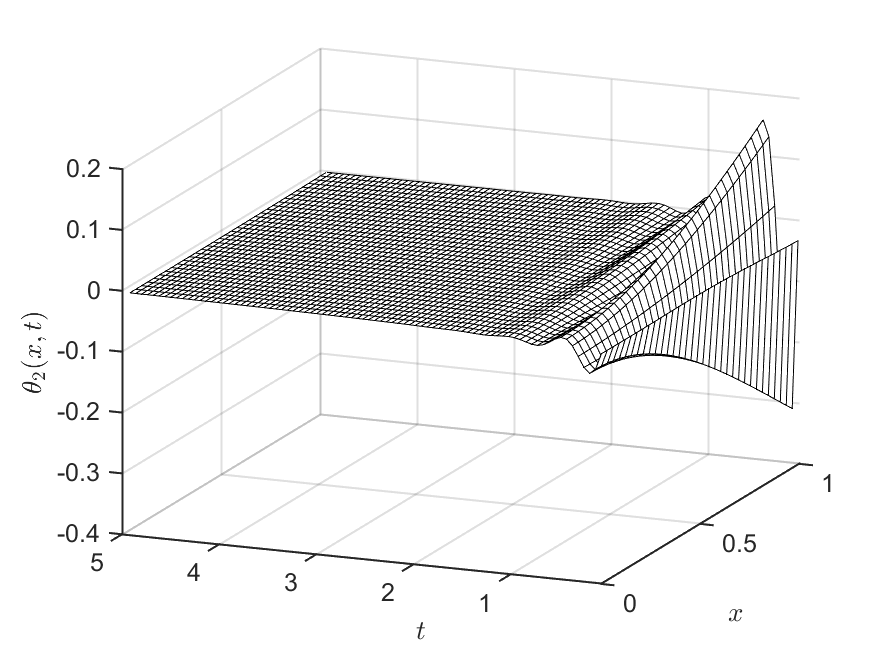}
		\caption{Evolution of closed-loop Timoshenko states $v_1(x,t),\theta_1(x,t), v_2(x,t),\theta_2(x,t)$ (from left to right).}
		\label{Fig-closed}
	\end{centering}
\end{figure*} 

\section{Numerical simulation}\label{num-simulation}
To illustrate the stabilization result with a numerical example,
we consider a particular case of multilayer Timoshenko beam with two layers ($N=2$) and with $\beta_1=1,\beta_2=2, \eta_1=\eta_2=1, \zeta_1=1, \zeta_2=2, h^1_1=1, h^1_2=1, k^1_T=1,k^1_N=1, \alpha_1=1, \alpha_2=1, \xi_1=\xi_3=-1, \xi_2=\xi_4=1$, where the transport speeds verify $\frac{\sqrt{\eta_1}}{\sqrt{\beta_1}}=\frac{\sqrt{\alpha_1}}{\sqrt{\zeta_1}}>\frac{\sqrt{\eta_2}}{\sqrt{\beta_2}}=\frac{\sqrt{\alpha_2}}{\sqrt{\zeta_2}}$. Thus, there are two ``blocks'' with the same transport speeds. The value of $\Phi(0)$ can be set arbitrarily to specifiy the decay rate and is chosen as $\Phi_{11}(0)=-11, \Phi_{12}(0)=1, \Phi_{22}(0)=-5, \Phi_{33}(0)=-11.4142, \Phi_{34}(0)=1, \Phi_{44}(0)=5$ and other elements in $\Phi(0)$ are set as zero, which leads to  the decay rate $C_2=(c'+1)/2-1=c-1=4$. For computing the 48 highly coupled kernel equations with triangle domains, a power series method \cite{Vazquez2024} is applied in which each kernel function is expressed by a polynomial of tenth order with undetermined coefficients that can be obtained in terms of equations and boundary conditions by the software Mathematica. The steps for the power series method are shown in the Fig. \ref{Fig-power}. Firstly, calculating the number $n_l$ of different characteristic lines whose slopes are less than 1 and greater than 0, and determining the number of area splits $n_l+1$. Secondly, defining the kernel functions as power series containing undetermined coefficients in each area split. Thirdly, inserting the power series expressions of kernel functions into the kernel equations and boundaries, comparing the coefficients of the same power terms and establishing the linear equations containing only undetermined parameters. After that, solving the linear equations regarding undetermined coefficients and obtaining the solutions of piecewise kernel functions. Finally, piecewise kernels are combined to obtain the final solutions. According to the characteristic lines of the kernel equations $K_{13},K_{14}, K_{23},K_{24}$ are piecewise functions and each of them needs two boundary conditions (one is at $y=0$ and another one is at $y=x$). Due to the complex coupling relationship between kernel functions, it is necessary to solve the remaining kernel functions as piecewise functions, and the continuities near the partition line can be used as the boundary conditions. Each kernel function has ten undetermined coefficients and can be determined by the rule that the coefficient of the term with the same order on both sides of the equal sign is also the same. Specifically, we start by solving $K_{1j}, L_{1j},\Phi_{1j}, j=1,2,3,4$ since they are coupled with each other and independent of other kernel functions.The kernel equations are solved taking into account that $K_{13},K_{14}$ are all discontinuous due to the fact that (\ref{ker_4}) and (\ref{ker_5}) needs to be simultaneously verified (but not $K_{12}$ due to it belonging to the same ``block'' as $K_{11}$). Thus, they possess a ``line of discontinuity'' along which they should be split in two analytic parts by dividing the triangular domain $\mathcal T$ into two parts which means we equivalently solve 2x4+2x4+4=20 coupled kernel functions. Next, $K_{2j}, L_{2j},\Phi_{2j}, j=1,2,3,4$ are solved. Since $K_{23},K_{24}$ are also discontinuous and they are coupled with $K_{1j}, L_{1j},\Phi_{1j}, j=1,2,3,4$, we need to divide the areas to two parts which means we equivalently solve 2x4+2x4+4=20 coupled kernel functions. This procedure is followed until all the kernels are found. We show the solutions of the gain kernels $K_{ij}(1,y),L_{ij}(1,y), 1\le i\le4,1\le j\le4$ in Fig. \ref{Fig-kernels}. 
The boundary values $\Phi$ appearing in the control law are 
$$\Phi(1)\hspace{-3pt}=\hspace{-5pt} \left[ \begin{array}{cccc} 
-11.4603 &  3.21899 &-4.64424& -0.62423\\ -7.07472 & 3.023 & 
-2.26582 &1.53245\\
-3.71605 & -1.41335 &-10.2616 &  -2.66744\\
	-3.968 & 6.03106 &-22.4318 & 30.2566
\end{array}\right]
$${
The numerical simulations reveal several important aspects of the system behavior. The open-loop response shown in Fig. 4 demonstrates the strong destabilizing effect of the boundary conditions, manifesting in exponential divergence across all states with pronounced spatial variations. Most notably, the second layer exhibits stronger divergence than the first, illustrating the complex coupling effects between layers. When the proposed controller is applied, Fig. 5 shows not only the achievement of exponential stability but also the preservation of the coupling structure during the stabilization process - the states maintain their relative amplitude relationships while converging uniformly to zero. The spatial variations in the response are effectively managed by the controller, which is particularly evident in the smooth evolution of the angular displacements $\theta_1$ and $\theta_2$. The simulation results thus validate both the theoretical stability guarantees and the controller's ability to handle the complex spatial-temporal dynamics inherent to multi-layer Timoshenko systems.}


\section{Concluding remarks}\label{sec-conclusion}
{This work presented boundary control of a $N$-layer Timoshenko composite beam with anti-damping and anti-stiffness at the uncontrolled boundary by extending the backstepping hyperbolic design to isotachic PIDE-ODE systems described in blocks. The $N$-layer Timoshenko beam can be written in such a way, and we are able to impose arbitrarily fast decay by applying the block PDE backstepping design. The numerical simulations demonstrate this conclusively: while the open-loop system exhibits rapid divergence with amplitudes growing by three orders of magnitude within the simulation window, the closed-loop response shows exponential convergence to zero with the theoretically predicted rate withing a few seconds. This performance was achieved through a complex control design involving 48 coupled kernel equations arranged in triangular domains, whose numerical solution required careful implementation of a tenth-order power series method. The controller here employs independent actuation of the boundary forces and torques at the boundary of each of the layers of the beam, as a consequence of the possible instability at the tip and of the rapid stabilization objective, not achievable with a single/common force and/or torque. Independent multiple actuations are also required for rapidly stabilizing controllers for  multi-layer fluids \cite{diagne2017control} and multi-lane traffic \cite{yu2021output}. It is only when the subsystems occupy the same physical space, as in multi-class traffic \cite{burkhardt2021stop}, that a single actuator suffices. As a possible future extension, beyond observer design based in the same principles, the ideas could be used to improve the design of the target system. The problem of designing observers and output feedback controllers will be considered in future work, with the goal of achieving arbitrarily rapid stabilization through output feedback control.}

\section*{Acknowledgments}

Vazquez was supported by grant TED2021-132099B-C33 funded by MICIU/ AEI/ 10.13039 /501100011033 and by ``European Union NextGenerationEU/PRTR.'' Chen and Qiao were supported by the National Key Research and Development Program of China under Grant 2021ZD0112301, the National Natural Science Foundation of China under Grant 62403018. Krstic was supported by AFOSR   FA9550-23-1-0535. 

\section*{Data Availability Statement}
The data that support the findings of this study are available on request from the corresponding author. The data are not publicly available due to privacy or ethical restrictions.
\bibliography{reference}

\begin{thebibliography}{10}
\providecommand \doibase [0]{http://dx.doi.org/}%

\bibitem{Ferreira2003}
Ferreira A. Thick composite beam analysis using a global meshless approximation
  based on radial basis functions. {\it Mech. Adv. Mater. Struct..}
  2003\string;10(3)\string:271--284.

\bibitem{Attwood2016}
Attwood J, Russell B, Wadley H, Deshpande V. Mechanisms of the penetration of
  ultra-high molecular weight polyethylene composite beams. {\it Int. J. Impact
  Eng..} 2016\string;93\string:153--165.

\bibitem{Barbero1991}
Barbero EJ, Fu SH, Raftoyiannis I. Ultimate bending strength of composite
  beams. {\it Journal of Materials in Civil Engineering.}
  1991\string;3(4)\string:292--306.

\bibitem{Waisman2002}
Waisman H, Abramovich H. Active stiffening of laminated composite beams using
  piezoelectric actuators. {\it Composite structures.}
  2002\string;58(1)\string:109--120.

\bibitem{mattioni2023lyapunov}
Mattioni A, Wu Y, Le~Gorrec Y. A Lyapunov approach for the exponential
  stability of a damped Timoshenko beam. {\it IEEE transactions on automatic
  control.} 2023.

\bibitem{zhao2020finite}
Zhao Z, Liu Z. Finite-time convergence disturbance rejection control for a
  flexible Timoshenko manipulator. {\it IEEE/CAA Journal of Automatica Sinica.}
  2020\string;8(1)\string:157--168.

\bibitem{guo2021robust}
Guo BZ, Meng T. Robust output regulation for Timoshenko beam equation with two
  inputs and two outputs. {\it International Journal of Robust and Nonlinear
  Control.} 2021\string;31(4)\string:1245--1269.

\bibitem{homaeinezhad2023feedback}
Homaeinezhad MR, Abbasi~Gavari M. Feedback control of actuation-constrained
  moving structure carrying Timoshenko beam. {\it International Journal of
  Robust and Nonlinear Control.} 2023\string;33(3)\string:1785--1806.

\bibitem{meng2023reduced}
Meng T, Wang J, Fu Q, He X. Reduced-order observer based output feedback
  control for a Timoshenko beam. {\it International Journal of Control.}
  2023\string:1--17.

\bibitem{wang2023saturated}
Wang S, Han ZJ, Zhao ZX, Cao XG. Saturated boundary feedback stabilization for
  a non-uniform Timoshenko beam with disturbances and uncertainties. {\it IEEE
  Transactions on Automatic Control.} 2023.

\bibitem{redaud2022domain}
Redaud J, Auriol J, Le~Gorrec Y. In-domain damping assignment of a
  Timoshenko-beam using state feedback boundary control. In: IEEE.
  2022\string:5405--5410.

\bibitem{Lo2015}
Lo A, Tatar NE. Stabilization of laminated beams with interfacial slip. {\it
  Electronic Journal of Differential Equations.} 2015\string;2015\string:129.

\bibitem{Mustafa2018}
Mustafa MI. Laminated {T}imoshenko beams with viscoelastic damping. {\it J.
  Math. Anal. Appl..} 2018\string;466(1)\string:619--641.

\bibitem{Apalara2020}
Apalara TA, Raposo CA, Nonato CA. Exponential stability for laminated beams
  with a frictional damping. {\it Archiv der Mathematik.}
  2020\string;114\string:471--480.

\bibitem{Apalara2020_1}
Apalara TA, Nass AM, Al~Sulaimani H. On a laminated {T}imoshenko beam with
  nonlinear structural damping. {\it Math. Comput. Appl..}
  2020\string;25(2)\string:35.

\bibitem{Guesmia2021}
Guesmia A, Mu{\~n}oz~Rivera JE, Sep{\'u}lveda~Cort{\'e}s MA, Vera~Villagr{\'a}n
  O. Laminated {T}imoshenko beams with interfacial slip and infinite memories.
  {\it Math. Methods Appl. Sci..} 2022\string;45(8)\string:4408--4427.

\bibitem{Cao2007}
Cao XG, Liu DY, Xu GQ. Easy test for stability of laminated beams with
  structural damping and boundary feedback controls. {\it Journal of Dynamical
  and Control Systems.} 2007\string;13(3)\string:313--336.

\bibitem{Tatar2015}
Tatar NE. Stabilization of a laminated beam with interfacial slip by boundary
  controls. {\it Boundary Value Problems.} 2015\string;2015(1)\string:1--11.

\bibitem{Alves2020}
Alves M, Monteiro R. Exponential stability of laminated {T}imoshenko beams with
  boundary/internal controls. {\it Journal of Mathematical Analysis and
  Applications.} 2020\string;482(1)\string:123516.

\bibitem{Feng2018}
Feng B. Well-posedness and exponential decay for laminated {T}imoshenko beams
  with time delays and boundary feedbacks. {\it Math. Methods Appl. Sci..}
  2018\string;41(3)\string:1162--1174.

\bibitem{krstic2008control}
Krstic M, Guo BZ, Balogh A, Smyshlyaev A. Control of a tip-force destabilized
  shear beam by observer-based boundary feedback. {\it SIAM Journal on Control
  and Optimization.} 2008\string;47(2)\string:553--574.

\bibitem{smyshlyaev2009arbitrary}
Smyshlyaev A, Guo BZ, Krstic M. Arbitrary decay rate for {E}uler-{B}ernoulli
  beam by backstepping boundary feedback. {\it IEEE Transactions on Automatic
  Control.} 2009\string;54(5)\string:1134--1140.

\bibitem{chen2022}
Chen G, Vazquez R, Krstic M. Rapid Stabilization of Timoshenko Beam by PDE
  Backstepping. {\it IEEE Transactions on Automatic Control.}
  2024\string;69(2)\string:1141-1148.

\bibitem{chen2023}
Chen G, Vazquez R, Krstic M. Backstepping-based Rapid Stabilization of
  Two-layer Timoshenko Composite Beams. {\it IFAC-PapersOnLine.}
  2023\string;56(2)\string:8159-8164.

\bibitem{krstic3}
Krstic M, Smyshlyaev A. Backstepping boundary control for first-order
  hyperbolic {PDE}s and application to systems with actuator and sensor delays.
  {\it Systems \& Control Letters.} 2008\string;57\string:750-758.

\bibitem{vazquez-linear}
Vazquez R, Krstic M, Coron JM. Backstepping Boundary Stabilization and State
  Estimation of a 2x2 Linear Hyperbolic System. {\it in Proceedings of the 50th
  CDC.} 2011.

\bibitem{vazquez-nonlinear}
Coron JM, Vazquez R, Krstic M, Bastin G. Local Exponential {$H^2$}
  Stabilization of a $2 \times 2$ Quasilinear Hyperbolic System using
  Backstepping. {\it SIAM J. Control Optim..} 2013\string;51\string:2005--2035.

\bibitem{florent}
Di~Meglio F, Vazquez R, Krstic M. Stabilization of a System of $n+1$ Coupled
  First-Order Hyperbolic Linear {P}{D}{E}s With a Single Boundary Input. {\it
  IEEE Transactions on Automatic Control.} 2013\string;58(12)\string:3097-3111.

\bibitem{hu}
Hu L, Di~Meglio F, Vazquez R, Krstic M. Control of homodirectional and general
  heterodirectional linear coupled hyperbolic {PDE}s. {\it IEEE Transactions on
  Automatic Control.} 2015\string;61(11)\string:3301--3314.

\bibitem{c6}
Hu L, Vazquez R, Meglio FD, Krstic M. Boundary exponential stabilization of
  1-dimensional inhomogeneous quasi-linear hyperbolic systems. {\it SIAM
  Journal on Control and Optimization.} 2019\string;57(2)\string:963--998.

\bibitem{Auriol2016}
Auriol J, Di~Meglio F. Minimum time control of heterodirectional linear coupled
  hyperbolic {PDE}s. {\it Automatica.} 2016\string;71\string:300--307.

\bibitem{espitia2021predictor}
Espitia N, Perruquetti W. Predictor-feedback prescribed-time stabilization of
  LTI systems with input delay. {\it IEEE Transactions on Automatic Control.}
  2021\string;67(6)\string:2784--2799.

\bibitem{c1}
Di~Meglio F, Argomedo FB, Hu L, Krstic M. Stabilization of coupled linear
  heterodirectional hyperbolic {PDE}--{ODE} systems. {\it Automatica.}
  2018\string;87\string:281--289.

\bibitem{deutscher2019output}
Deutscher J, Gehring N, Kern R. Output feedback control of general linear
  heterodirectional hyperbolic PDE-ODE systems with spatially-varying
  coefficients. {\it International Journal of Control.}
  2019\string;92(10)\string:2274--2290.

\bibitem{gabriel2023robust}
Gabriel J, Deutscher J. Robust Cooperative Output Regulation for Networks of
  Hyperbolic PIDE–ODE Systems. {\it IEEE Transactions on Automatic Control.}
  2024\string;69(2)\string:888-903.

\bibitem{humaidi2018design}
Humaidi A, Hameed M, Hameed A. Design of block-backstepping controller to ball
  and arc system based on zero dynamic theory. {\it Journal of Engineering
  Science and Technology.} 2018\string;13(7)\string:2084--2105.

\bibitem{humaidi2017design}
Humaidi AJ, Hameed MR. Design and performance investigation of
  block-backstepping algorithms for ball and arc system. In: IEEE.
  2017\string:325--332.

\bibitem{humaidi2021block}
Humaidi AJ, Hameed MR, Ajel AR, Hameed AH, Al-Qassar AA, Ibraheem IK. Block
  backstepping control design of two-wheeled inverted pendulum via zero dynamic
  analysis. In: IEEE.  2021\string:87--92.

\bibitem{Lenci2013}
Lenci S, Rega G. A Limit Model for the Linear Dynamics of a Two-Layer Beam.
  {\it International Design Engineering Technical Conferences and Computers and
  Information in Engineering Conference.} 2013\string;55997\string:V008T13A067.

\bibitem{Wonham1967}
Wonham W. On pole assignment in multi-input controllable linear systems. {\it
  IEEE Trans. Aut. Contr..} 1967\string;12(6)\string:660--665.

\bibitem{c2}
Su L, Wang JM, Krstic M. Boundary feedback stabilization of a class of coupled
  hyperbolic equations with nonlocal terms. {\it IEEE Transactions on Automatic
  Control.} 2017\string;63(8)\string:2633--2640.

\bibitem{Vazquez2024}
Vazquez R, Chen G, Qiao J, Krstic M. The Power Series Method to Compute
  Backstepping Kernel Gains: Theory and Practice.  2023\string:8162-8169.

\bibitem{diagne2017control}
Diagne M, Tang SX, Diagne A, Krstic M. Control of shallow waves of two unmixed
  fluids by backstepping. {\it Annual Reviews in Control.}
  2017\string;44\string:211--225.

\bibitem{yu2021output}
Yu H, Krstic M. Output feedback control of two-lane traffic congestion. {\it
  Automatica.} 2021\string;125\string:109379.

\bibitem{burkhardt2021stop}
Burkhardt M, Yu H, Krstic M. Stop-and-go suppression in two-class congested
  traffic. {\it Automatica.} 2021\string;125\string:109381.

\end{thebibliography}

\end{document}